\documentclass[reqno]{amsart}
\usepackage{amsmath,amsxtra,amssymb,latexsym, amscd,amsthm}
\usepackage[mathscr]{euscript}
\usepackage{mathrsfs}
\usepackage{color}
\usepackage[active]{srcltx}
\usepackage{mathtools,stmaryrd}
\usepackage{enumerate}
\usepackage{hyperref}
\usepackage{framed}
\usepackage[T1]{fontenc}
\usepackage{subfigure}
\usepackage{multirow}
\usepackage{booktabs}
\usepackage{tikz}
\usepackage{longtable}
\usepackage{tikz-qtree}
\usepackage[edges]{forest}
\usetikzlibrary{trees} 
\usepackage{stmaryrd}
\SetSymbolFont{stmry}{bold}{U}{stmry}{m}{n}
\SetSymbolFont{stmry}{bold}{U}{stmry}{b}{n}
\usepackage{graphicx}
\usepackage{subcaption}

\newcommand{\RR}{{\mathbb R}}
\newcommand{\CC}{{\mathbb C}}

\newcommand{\bS}{{\mathbb S}}

\newcommand{\N}{\mathbb{N}}

\newcommand{\K}{\mathcal{K}}
\newcommand{\MQ}{\mathcal{Q}}

\newcommand{\bG}{\mathbf{G}}

\newcommand{\bY}{\mathbf{Y}}

\def\bx{{\boldsymbol{x}}}

\def\by{{\boldsymbol{y}}}

\def\bz{{\boldsymbol{z}}}

\def\bv{{\boldsymbol{v}}}
\def\bu{{\boldsymbol{u}}}

\newcommand{\supp}{\mathrm{supp}}

\newcommand{\lram}[3]{\left\langle #1, #2\right\rangle_{#3}}

\def\bfm{{\boldsymbol{m}}}

\def\bu{{\boldsymbol{u}}}

\def\bv{{\boldsymbol{v}}}

\def\bgm{{\boldsymbol{\gamma}}}

\newcommand{\ud}{\mathrm{d}}

\newcommand{\rank}{\mbox{\upshape rank}}

\newcommand{\p}{\mathsf{p}}
\newcommand{\ch}{\mathsf{ch}}

\newcommand{\cssos}{\mbox{\upshape\tiny cs-sos}}
\newcommand{\csmom}{\mbox{\upshape\tiny cs-mom}}

\newcommand{\tr}[1]{\mbox{\upshape tr}\left(#1\right)}
\def\ba{{\boldsymbol{\alpha}}}
\def\bb{{\boldsymbol{\beta}}}

\newcommand{\cA}{\mathcal{A}}
\newcommand{\mA}{\mathcal{A}}

\newcommand{\mL}{\mathscr{L}}
\newcommand{\mI}{\mathcal{I}}

\newcommand{\mJ}{\mathcal{J}}
\newcommand{\wt}[1]{\widetilde{#1}}

\newcommand{\mG}{\mathcal{G}}

\newcommand{\mE}{\mathcal{E}}

\newcommand{\QM}{\mathcal{Q}}

\newcommand{\X}{\mathcal{X}}

\newtheorem{assumption}{Assumption}[section]

\newtheorem{theorem}{Theorem}[section]

\newtheorem{corollary}{Corollary}[section]
\newtheorem{lemma}{Lemma}[section]
\newtheorem{example}{Example}[section]
\newtheorem{definition}{Definition}[section]
\newtheorem{remark}{Remark}[section]
\newtheorem{proposition}{Proposition}[section]

\title{
A Correlative-Sparse Matrix Positivstellensatz under Pointwise Commutation with Applications to Polynomial Matrix Optimization
}

\author{Feng Guo}
\address[Feng Guo]{School of Mathematical Sciences, Dalian University of Technology, Dalian 116024, Liaoning Province, China.}
\email{fguo@dlut.edu.cn}

\thanks{The author was supported by the National Natural Science Foundation of China (grant no. 12471478).}

\date{\today}

\begin{document}

\begin{abstract}
The scalar Positivstellensatz adapted to correlative sparsity enables substantially smaller moment--SOS relaxations for sparse polynomial optimization. Its matrix counterpart, however, fails in general, even under the running intersection property and local Archimedean assumptions. In this paper, we establish a matrix Positivstellensatz under these assumptions and an additional pointwise commutation condition. Specifically, if a polynomial matrix admits a decomposition whose summands commute pairwise on the feasible set, then positive definiteness on this set implies membership in the sum of the local matrix quadratic modules. The proof uses central common lower envelopes to derive an integral inequality for positive matrix-valued clique measures with consistent separator marginals, without requiring these measures to admit a global extension. We apply this result to establish asymptotic convergence of the corresponding sparse moment--SOS hierarchy for polynomial matrix optimization. We also derive a criterion for finite convergence of moment relaxation and extraction of a global minimizer based on local flatness and consistency of the local representing matrix measures on overlapping variable blocks. Finally, we provide a sufficient construction of commuting decompositions through linear equations, 
using a reference polynomial matrix with a simple spectrum at some point.
\end{abstract}


\maketitle

\section{Introduction}
Positivstellens\"atze are fundamental results in real algebraic geometry,
providing algebraic certificates for polynomial positivity on semialgebraic sets
\cite{realAG,Laurent_sumsof,PPSOSMarshall,ScheidererSurvey}.
These results have widespread applications in optimization, algebraic geometry,
and control theory, where polynomial positivity plays a central role
\cite{SOCAGbook,HKL2020,Lasserre2015,Laurent_sumsof,Nie2023}.

Typically, Positivstellens\"atze characterize polynomials
$f(\bx)\in\RR[\bx]$ that are positive on a basic semialgebraic set
\begin{equation}\label{eq::set1}
\{\bx\in\RR^n \mid g_j(\bx)\ge 0,\ j=1,\ldots,s\},
\end{equation}
where $\RR[\bx]$ denotes the ring of real polynomials in
$\bx=(x_1,\ldots,x_n)$ and $g_1,\ldots,g_s\in\RR[\bx]$.
Many such results express positivity through sums of squares (SOS).
For instance, if the set \eqref{eq::set1} is compact,
Schm\"udgen \cite{Schmugen1991} showed that every polynomial $f$
positive on this set belongs to the preordering generated by
$g_1,\ldots,g_s$; that is, it admits an SOS-weighted representation
using products of these polynomials.
Under the Archimedean assumption, which is slightly stronger than compactness,
Putinar's Positivstellensatz \cite{Putinar1993} provides a quadratic
module representation that avoids such products.
Building on this result, Lasserre
\cite{LasserreGlobal2001,Lasserre09,Lasserre2015}
introduced the well-known hierarchy of moment--SOS relaxations
for polynomial optimization.

SOS-based characterizations of polynomial positivity extend to the matrix setting.
Using SOS polynomial matrices, Scherer and Hol \cite{SH2006}
developed, under an Archimedean assumption, a matrix analogue of
Putinar's Positivstellensatz for an $m\times m$ symmetric polynomial
matrix $F(\bx)$ that is positive definite (PD) on a set defined by
polynomial matrix inequalities (PMIs):
\begin{equation}\label{eq::K}
    \K:=\{\bx\in\RR^n \mid
    G_1(\bx)\succeq 0,\ldots,G_s(\bx)\succeq 0\},
\end{equation}
where each $G_j(\bx)$ is a $q_j\times q_j$ symmetric polynomial
matrix, $j=1,\ldots,s$.
Dinh et al. \cite{DHL2021} further generalized the classical
Schm\"udgen, Putinar--Vasilescu, and Dickinson--Povh
Positivstellens\"atze to the polynomial matrix setting.
For the projections of two PMI-defined semialgebraic sets,
Klep and Nie \cite{KN2020} provided a matrix Positivstellensatz
involving lifting polynomials to certify set containment.
These generalizations have a wide range of applications,
particularly in optimal control and systems theory
\cite{HL2006,HL2012,ichihara2009optimal,pozdyayev2014atomic,vanantwerp2000tutorial}.

Computing the SOS representations appearing in Positivstellens\"atze
can typically be formulated as semidefinite programs (SDPs).
Their computational cost depends strongly on the sizes and number
of positive semidefinite (PSD) matrix constraints
\cite{lasserre2009moments}.
For an $m\times m$ SOS polynomial matrix of degree at most $2d$
in $n$ variables, the dense Gram matrix has size
$m\binom{n+d}{d}$.
This size grows linearly with $m$ but rapidly with $n$ and $d$,
leading to substantial computational costs.

Structural information can reduce the cost of SOS verification
by shrinking monomial bases and decomposing PSD constraints
\cite{magron2023sparse,zheng2021chordal}.
In the scalar setting, the main approaches include correlative
sparsity \cite{waki2006sums}, term sparsity in commutative
\cite{wang3,wang2021tssos} and noncommutative
\cite{wang2021exploiting} settings, and symmetry reduction
\cite{gatermann2004symmetry}.
Newton polytopes eliminate unnecessary monomials in unconstrained
SOS representations \cite{Reznick_1978}, while term sparsity is
closely related to sign-symmetry reduction
\cite{wang2021tssos,lofberg2009pre}.
Equality constraints also enable algebraic reductions
\cite{parrilo2005exploiting} and ideal-sparsity decompositions.
These structures can be combined, as in CS-TSSOS \cite{wang2022cs}.
Term sparsity has also been applied to sums of rational functions
\cite{guo2024exploiting} and non-SOS nonnegativity certificates.
In the matrix setting, chordal matrix sparsity yields sparse SOS
representations \cite{zheng2023sum}, independently of monomial
sparsity within nonzero entries.
For scalar objectives with matrix constraints, existing work
exploits correlative sparsity \cite{KM2009,nie2024} and 
constraint matrix sparsity \cite{kim2011exploiting}.
%and term sparsity in frame topology optimization \cite{handaterm}.

In our previous work \cite{MWG2026}, we showed how to exploit
different types of sparsity---term sparsity, correlative sparsity,
and matrix sparsity---in polynomial matrices to derive sparse SDP
relaxations for polynomial matrix optimization problems.
For correlative sparsity, however, we provided a counterexample (\cite[Example~5.1]{MWG2026})
showing that, unlike the scalar case, the corresponding matrix
Positivstellensatz does not hold in general.
Specifically, suppose that the variables $\bx$ are grouped into
finitely many, possibly overlapping blocks, such that each monomial
of the objective matrix $F$ and each constraint polynomial
matrix $G_j$ depends only on the variables in a single block.
Even if $F$ is positive definite on $\K$, it need not belong to
the sum of the local quadratic modules generated by the constraints
assigned to each block, viewed in the ambient polynomial ring,
despite the running intersection property (RIP) and the local Archimedean assumption.
This raises the question of what additional conditions on
$(F,G_1,\ldots,G_s)$ ensure a correlative-sparse matrix
Positivstellensatz.

\vskip 5pt
\noindent{\itshape Contributions.}
The main contributions of this paper are as follows.

(i) We establish a matrix Positivstellensatz adapted to correlative sparsity under a pointwise commutation condition. Specifically, suppose that $F$ admits a decomposition into polynomial matrices, each depending only on the variables in one block, whose summands commute pairwise on $\K$. Under the RIP and the local Archimedean assumption, positive definiteness of $F$ on $\K$ implies that $F$ belongs to the sum of the local quadratic modules. The proof relies on an integral inequality for consistent matrix-valued clique measures, without requiring a global representing measure.

(ii) We apply this Positivstellensatz to polynomial matrix optimization and establish the asymptotic convergence of the moment--SOS hierarchy adapted to correlative sparsity under the same structural assumptions. We also provide a criterion for finite convergence of moment relaxation and minimizer extraction based on local flatness and consistency of the local representing matrix measures on overlapping variable blocks. This criterion does not require these measures to admit a global extension. Finally, using a reference polynomial matrix with a simple spectrum at some point, we give a sufficient condition for constructing a commuting decomposition of $F$ (if it exists) by solving linear equations.

\vskip 7pt
The rest of this paper is organized as follows.  
We recall some notation and preliminaries needed for this paper in Section \ref{sec:preliminaries}.
A correlative sparse Positivstellensatz is established in Section \ref{sec::csPos}.
We apply this result to polynomial matrix optimization in Section \ref{sec::pmo}.
Some conclusions are summarized in Section \ref{sec::cons}.

\section{Preliminaries}
\label{sec:preliminaries}
We collect some notation and basic preliminaries
which will be used in this paper. 
\subsection{Notation}
The $n$-dimensional real Euclidean vector space is $\RR^n$. 
Let $\N=\{0,1,2,\ldots\}$, and let $\N^n$ denote the set of
$n$-dimensional multi-indices.
%The set of natural numbers is $\N$, and the $n$-dimensional set of multi-indices is $\N^n$. 
%\[
%\mathbb N_d^n=\{\boldsymbol\alpha\in\mathbb N^n:|\boldsymbol\alpha|\le d\},
%\qquad
%\mathbb N_d^S=\{\boldsymbol\alpha\in\mathbb N_d^n:
%\alpha_j=0\text{ for all }j\notin S\}.
%\]
The symbol $[n]$ denotes the set $\{1,\ldots,n\}$. For a set $\mA$, its cardinality is denoted by $|\mA|$. The degree of a multi-index $\ba \in \N^n$ is $|\ba| \coloneqq \sum_{i \in [n]} \alpha_i$.
For $d\in\N$, let $\N_d^n=\{\ba\in\N^n\mid |\ba|\le d\}$. 
The associated monomial to an $n$-dimensional indeterminate $\bx\coloneqq(x_1,\ldots,x_n)$ and a multi-index $\ba \in \N^n$ is $\bx^\ba = \prod_{i=1}^n x_i^{\alpha_i}$. The set of polynomials with real-valued coefficients is $\RR[\bx]$. The \emph{support} of a polynomial $f \in \RR[\bx]$, denoted by $\supp(f)$, is the set of multi-indices $\ba$ such that $\bx^{\ba}$ appears in $f$.
The degree of a polynomial $f \in \RR[\bx]$ is $\deg(f) = \max_{\ba \in \supp(f)} |\ba|$.
The set of polynomials of degree at most $d$ is denoted by $\RR[\bx]_{d}$. 
The set of $m\times m$ complex-valued matrices is denoted by $M_m(\CC)$.
The transpose of a matrix $A\in M_m(\CC)$ is denoted by $A^{\intercal}$, the Hermitian transpose
of $A$ is denoted by $A^*$ and the entrywise complex conjugation of $A$ is denoted by $\overline{A}$. 
For $A\in M_m(\CC)$, $\Vert A\Vert$ denotes the spectral norm of $A$, $\tr{A}$ denotes
the trace of $A$. Denote by $\lambda_{\min}(A)$ (resp., $\lambda_{\max}(A)$) the smallest (resp., 
largest) eigenvalue of a Hermitian $A\in M_m(\CC)$.
The $m$-dimensional identity matrix is $I_m$. For two matrices $C$ and $D$, 
$C \otimes D$ denotes the Kronecker product of $C$ and $D$, and $C \bigoplus D$ denotes
the direct sum of $C$ and $D$. 
The set of $m\times m$ symmetric and real-valued matrices is $\bS^m$. The subset of $m\times m$ PSD matrices is $\bS^m_+$. Membership in the PSD set $Z \in \bS^m_+$ will also be denoted as $Z \succeq 0$. 
The set of polynomial matrices of dimension $m \times m$ is $\RR[\bx]^{m \times m}$. Given a polynomial matrix $Z(\bx)=[Z(\bx)_{ij}]\in\RR[\bx]^{p \times q}$, we can write it as $Z(\bx)=\sum_{\ba}Z_{\ba}\bx^\ba$, $Z_{\ba}\in\RR^{p \times q}$. The support of $Z(\bx)$, denoted by $\supp(Z)$, is the set of multi-indices $\ba$ such that $Z_\ba \neq 0$.
The degree of $Z(\bx)$ is $\deg{Z} = \max_{i \in [p], \ j \in [q]} \deg{Z(\bx)_{ij}}$.
The set of symmetric polynomial matrices of dimension $m$ is $\bS^m[\bx]$. The set of $m$-dimensional symmetric polynomial matrices of degree at most $d$ is denoted by $\bS^m[\bx]_{d}$.
A graph is denoted by $\mathcal G=(V,E)$, where $V$ is a finite vertex set and
$E\subseteq\{\{i,j\}:i,j\in V,\ i\ne j\}$ is the edge set.
Throughout this paper, all graphs are simple and undirected.
A tree, denoted by $\mathcal T=(V,E)$, is a connected graph containing no cycles.
When necessary, we use $V_{\mathcal G},E_{\mathcal G}$ and
$V_{\mathcal T},E_{\mathcal T}$ to distinguish their vertex and edge sets.

\subsection{Sum of Squares Polynomial Matrices}
A polynomial matrix $Z(\bx) \in \bS^m[\bx]$ is positive semidefinite (PSD) over $\RR^n$ if $Z(\bx) \succeq 0$ for all $\bx \in \RR^n$. 
Furthermore, we call $Z(\bx)$ an SOS matrix if there exists another polynomial matrix $R(\bx) \in \RR[\bx]^{r \times m}$ such that $Z(\bx) = R(\bx)^{\intercal} R(\bx)$.
The set of $m$-dimensional SOS matrices is denoted by $\Sigma^m[\bx]$. For $k\in\N$,
let $\Sigma^m[\bx]_{2k}=\Sigma^m[\bx]\cap\bS^m[\bx]_{2k}$.
Note that the set of SOS (scalar) polynomials $\Sigma[\bx]$ is equivalent to the set $\Sigma^1[\bx]$ (the $m=1$ case of SOS matrices). 

Checking if a polynomial matrix $Z(\bx) \in \bS^m[\bx]_{2d}$ is inside the set $\Sigma^m[\bx]$ can be accomplished by solving an SDP. Letting $\bfm_d(\bx)$ be the $\binom{n+d}{d}$-vector of monomials in $\bx$ up to degree $d$, the matrix $Z(\bx)$ is an SOS matrix if there exists a PSD matrix $Q \in \bS_+^{m \binom{n+d}{d}}$ such that
{\small\begin{equation}
    Z(\bx) = \left(\bfm_d(\bx) \otimes I_m\right)^{\intercal} Q (\bfm_d(\bx) \otimes I_m). \label{eq:sos_pmi_uncons}
\end{equation}}
Checking the existence of $Q$ in \eqref{eq:sos_pmi_uncons} is an SDP with respect to the PSD matrix constraint $Q \in \bS_+^{m \binom{n+d}{d}}$ and the affine equality (coefficient matching) constraints in \eqref{eq:sos_pmi_uncons}.

Let $\bG:=\{G_j(\bx)\}_{j=1}^{s}$ be a set of symmetric polynomial matrices with $G_j(\bx) \in \bS^{q_j}[\bx]$ which defines a basic semialgebraic set as
\[
    \K \coloneqq \{\bx \in \RR^n \mid G_j(\bx) \succeq 0, \quad \forall j \in [s]\}.
\]
A polynomial matrix $Z(\bx) \in \bS^m[\bx]$ satisfies PMI over the constrained region $\K$ if
{\small\begin{align}
    \forall \bx \in \K: \qquad Z(\bx) \succeq 0. \label{eq:pmi_cons}
\end{align}}
The following exposition of SOS verification of PMI originates from \cite{SH2006}.

Let $C$ be a $mq\times mq$ block matrix  described as
{\small\begin{align}
    C = \begin{bmatrix}
        C_{11} & C_{12} & \hdots & C_{1m} \\
        C_{21} & C_{22} & \hdots & C_{2m} \\
        \vdots & \vdots & \ddots & \vdots \\
        C_{m1} & C_{m2} & \hdots & C_{mm}
    \end{bmatrix} \label{eq:c_block}
\end{align}}
in which each block $C_{ij}$ is of size $q\times q$.
The $m$-product between $C$ and a $q\times q$ matrix $D$ is defined by
\[
\begin{aligned}
    \langle C, D\rangle_m\coloneqq\begin{bmatrix}
        \text{tr}(C_{11}^{\intercal}D) & \text{tr}(C_{12}^{\intercal}D) & \hdots & \text{tr}(C_{1m}^{\intercal}D) \\
        \text{tr}(C_{21}^{\intercal}D) & \text{tr}(C_{22}^{\intercal}D) & \hdots & \text{tr}(C_{2m}^{\intercal}D) \\
        \vdots & \vdots & \ddots & \vdots \\
        \text{tr}(C_{m1}^{\intercal}D) & \text{tr}(C_{m2}^{\intercal}D) & \hdots & \text{tr}(C_{mm}^{\intercal}D)
    \end{bmatrix}.
\end{aligned}
\]
The \emph{matrix quadratic module} $\QM^m(\bG)$ generated by $\bG$ is defined as the cone
\[
	\QM^m(\bG)\coloneqq\left\{S_0(\bx)+\sum_{j=1}^s\lram{S_j(\bx)}{G_j(\bx)}{m}\ 
 \middle\vert \
S_0\in\Sigma^m[\bx],\ S_j\in\Sigma^{mq_j}[\bx], \ j\in[s]
\right\}.
\]
The $k$-th truncated matrix quadratic module $\QM^m_k(\bG)$ is a subcone of $\QM^m(\bG)$ defined as 
\[
\QM^m_k(\bG)\coloneqq\left\{S_0(\bx)+\sum_{j=1}^s\lram{S_j(\bx)}{G_j(\bx)}{m}\ 
 \middle\vert \
    S_0\in\Sigma^m[\bx]_{2k}, S_j\in\Sigma^{mq_j}[\bx]_{2(k-d_j)}, \ 
d_j=\left\lceil\frac{\deg G_j}{2}\right\rceil,\ j\in[s]
\right\}.
\]
\begin{assumption}\label{assump_arch}
{\rm
%The quadratic module 
$\QM^m(\bG)$ is {\itshape Archimedean}, i.e.,
there is $C>0$ such that $C - \Vert \bx\Vert^2\in\QM^1(\bG)$.}
\end{assumption}
\begin{theorem}{\upshape (Scherer-Hol's Positivstellensatz \cite{SH2006})}\label{th::psatz}
Let Assumption \ref{assump_arch} hold
and $F(\bx)\in\bS^m[\bx]$ be positive definite on $\K$.
%$\{\bx\in\RR^n \mid H_j(\bx)\succeq 0, \ j\in[t]\}$. 
Then $F(\bx)\in\QM^m(\bG)$.
\end{theorem}

\subsection{Matrix-Valued Measures and Moment Matrices}
We recall some basics about the theory of matrix-valued measures and the
associated moment matrices, which will be used to establish conditions for detecting global optimality and extracting optimal solutions from the Moment-SOS relaxations of PMO problems. 

For the set $\K$ in \eqref{eq::K}, let
$\mathscr{B}(\K)$ be the smallest $\sigma$-algebra generated from the open subsets of $\K$
and $\mathfrak{m}(\K)$ the set of all finite signed Borel measures on $\K$. A measure
$\phi\in\mathfrak{m}(\K)$ is \emph{positive} if $\phi(\cA)\ge 0$ for all $\cA\in
\mathscr{B}(\K)$. Denote by $\mathfrak{m}_+(\K)$ the set of all finite positive Borel
measures on $\K$. For
$\phi\in\mathfrak{m}(\K)$, its support is defined as the support of its
total variation measure $|\phi|$. 
We write
\[
L^1(\K,|\phi|)
:=
\left\{
f:\K\to\RR \ \middle|\
f\text{ is Borel measurable and }
\int_{\K}|f|\,\ud|\phi|<\infty
\right\}.
\]

Let $\phi_{ij}\in\mathfrak{m}(\K)$, $i,j= 1,\ldots, m$. The $m\times m$
matrix-valued measure $\Phi$ on $\K$ is defined as the matrix-valued function
$\Phi\colon \mathscr{B}(\K) \to \RR^{m\times m}$ with
\[\Phi(\cA)\coloneqq[\phi_{ij}(\cA)]\in \RR^{m\times m}, \quad \forall \cA\in
\mathscr{B}(\K). 
\]
A matrix-valued measure $\Phi=[\phi_{ij}]$ is called symmetric if
$\phi_{ij}=\phi_{ji}$ for all $i,j$. It is called positive semidefinite
if it is symmetric and
\[
\Phi(\mathcal A)\succeq 0
\qquad\text{for every }\mathcal A\in \mathscr{B}(\K).
\]
The set $\supp(\Phi)\coloneqq\bigcup_{i,j=1}^m
\supp(\phi_{ij})$ is called the support of $\Phi$. 
%The matrix-valued integral of $Z\in\RR[\bx]$ with respect to the measure $\Phi$ is
%defined by 
%\[\int_{\K} Z(\bx) \ud \Phi(\bx)\coloneqq\left[\int_{\K}Z(\bx)\ud
%\phi_{ij}(\bx)\right]_{i,j=1,\ldots,p}\in\RR^{p\times p}.\]
We denote by $\mathfrak{M}^m(\K)$ (resp. $\mathfrak{M}^m_+(\K)$) the set of all $m\times m$ (resp. PSD) symmetric matrix-valued measures on $\K$. 
A \emph{finitely atomic PSD} matrix-valued measure $\Phi\in\mathfrak{M}^m_+(\K)$ is a matrix-valued measure of form $\Phi=\sum_{i=1}^t W_i \delta_{\bx^{(i)}}$ 
where $W_i\in\bS_+^m$, $\bx^{(i)}$'s are distinct points in $\K$, and $\delta_{\bx^{(i)}}$ denotes the Dirac measure centered at $\bx^{(i)}$, $i=1,\ldots,t$.

For a finite symmetric positive semidefinite matrix-valued measure
$\Phi=[\phi_{ij}]\in\mathfrak{M}_+^m(\K)$ and a Borel measurable
matrix-valued function $A=[a_{ij}]:\K\to\bS^m$ satisfying
$a_{ij}\in L^1(\K,|\phi_{ij}|)$, define the duality pairing
\begin{equation}\label{eq::Phi}
\langle A,\Phi\rangle
:=
\int_{\K}\operatorname{tr}\bigl(A(\bx)\,\ud\Phi(\bx)\bigr)
:=
\sum_{i,j=1}^m
\int_{\K}a_{ij}(\bx)\,\ud\phi_{ij}(\bx).
\end{equation}
Then, it is clear that if $A(\bx)\succeq 0$ on $\K$ and $\Phi\in\mathfrak{M}^m_+(\K)$, then 
$\langle A,\Phi\rangle\ge 0$.

Let $\bY=(Y_{\ba})_{\ba\in\N^n}$ be a multi-indexed sequence of symmetric matrices in $\bS^m$.
We define a linear functional $\mL_{\bY}  \colon \bS^m[\bx] \to \RR$ in the following way:
\[
\mL_{\bY}(F)\coloneqq\sum_{\ba\in\supp(F)}\tr{F_{\ba}Y_{\ba}},\quad \forall 
F(\bx)=\sum_{\ba\in\supp(F)} F_{\ba}\bx^{\ba}\in\bS^m[\bx].
\]
We call $\mL_{\bY}$ the \emph{Riesz functional} associated to the sequence $\bY$.
We say the sequence $\bY$ has a representing matrix-valued measure 
$\Phi=[\phi_{ij}]\in\mathfrak{M}^m_+(\K)$ if
\begin{equation}\label{eq::S}
\supp(\Phi)\subseteq \K\quad\text{and}\quad 
Y_{\ba}=\int_{\K} \bx^{\ba}\ud \Phi(\bx)\coloneqq\left[\int_{\K}\bx^{\ba}\ud
\phi_{ij}(\bx)\right]_{i,j\in[m]},\ \forall \ba\in\N^n.
\end{equation}
%The measure $\Phi\in\mathfrak{M}^m_+(\K)$ satisfying \eqref{eq::S} is called a \emph{representing measure} of $\bY$.
%For a given sequence $\bS=(S_{\ba})_{\ba\in\N^n}\subseteq\mS^{p}$, 
%Clearly, $\bS$ is a matrix-valued $\K$-moment sequence if and only if $\mL_{\bS}$ is a tracial $\K$-moment functional.

\begin{definition}\label{def::mm}
Given a sequence $\bY=(Y_{\ba})_{\ba\in\N^n}\subseteq\bS^m$,
the associated \emph{moment matrix} $M(\bY)$ is the block matrix whose block row and block column are 
indexed by $\N^n$ and the $(\ba, \bb)$-th block entry is $Y_{\ba+\bb}$ for all $\ba, \bb\in\N^n$. For
$G\in\bS^q[\bx]$,
the \emph{localizing matrix} $M(G\bY)$ associated to $\bY$ and $G$ is the block matrix whose block row and block column are indexed by $\N^n$ and 
the $(\ba, \bb)$-th block entry is 
$\sum_{\bgm\in\supp(G)}Y_{\ba+\bb+\bgm}\otimes G_{\bgm}$ for all $\ba, \bb\in\N^n$. For $d\in\N$, the $d$-th order moment matrix $M_d(\bY)$ (resp. localizing matrix
$M_d(G\bY)$) is the submatrix of $M(\bY)$ (resp. $M(G\bY)$) whose block row and block column are both indexed by $\N^n_d$.
\end{definition}

If a sequence $\bY\subseteq\bS^m$ has a representing matrix-valued measure 
$\Phi=[\phi_{ij}]\in\mathfrak{M}^m_+(\K)$, it is easy to check that
$M(\bY)\succeq 0$ and $M(G_j\bY)\succeq 0$ for all $j\in[s]$. Under the Archimedean assumption,
these conditions are also sufficient.

\begin{theorem}\label{th::MM}{\upshape \cite[Theorems 5 and 6]{CIMPRIC2013}}
Let Assumption \ref{assump_arch} hold. Given a sequence 
$\bY=(Y_{\ba})_{\ba\in\N^n}\subseteq\bS^{m}$, $\bY$ has a representing matrix-valued measure 
$\Phi=[\phi_{ij}]\in\mathfrak{M}^m_+(\K)$ if and only if $M(\bY)\succeq 0$ and $M(G_j\bY)\succeq 0$ for all $j\in[s]$.
\end{theorem}
We can identify each linear functional
$\mL:\bS^m[\bx]\to\RR$ with the unique sequence
$\bY=(Y_\ba)_{\ba\in\N^n}$ of matrices in $\bS^m$ satisfying
\[
\mL(\bx^\ba B)=\tr{B Y_\ba}
\qquad\text{for all } B\in\bS^m \text{ and } \ba\in\N^n.
\]
A measure $\Phi\in\mathfrak{M}^m_+(\K)$ represents $\mL$ if it
represents the associated sequence $\bY$, or equivalently, if
\[
\mL(F)=\langle F,\Phi\rangle
\qquad\text{for every } F\in\bS^m[\bx].
\]
Under this identification, the matrix inequalities in
Theorem~\ref{th::MM} are equivalent to the condition
\[
\mL(H)\ge 0 \qquad\text{for every } H\in\MQ^m(\bG).
\]
Thus, under Assumption~\ref{assump_arch}, this condition is necessary
and sufficient for $\mL$ to admit a representing measure in
$\mathfrak{M}^m_+(\K)$.

Let 
\[
d_i=\lceil\deg{G_i}/2\rceil,\ i\in[s],\quad
d_\bG\coloneqq\max\{1, d_1,\ldots,d_s\}.
\]
The following theorem gives a flat extension criterion for truncated
matrix-valued sequences.
\begin{theorem}[\cite{GW2025}, Theorem 5]\label{th::FEC}
Given a truncated sequence $\bY=(Y_{\ba})_{\ba\in\N^n_{2k}}\subseteq\bS^{m}$, 
the following statements are equivalent:
\begin{enumerate}[\upshape (i)]
\item $\bY$ admits an atomic representing measure $\Phi=\sum_{i=1}^r W_i\delta_{\bx^{(i)}}$
with $W_i\in\bS_+^{m}$, $\bx^{(i)}\in\K$ and $\sum_{i=1}^r\rank(W_i)=\rank(M_k(\bY))$;
\item $M_k(\bY)\succeq 0$ and $\bY$ admits an extension $\wt{\bY}=(\wt{Y}_{\ba})_{\ba\in\N^n_{2(k+d_{\bG})}}$
such that $M_{k+d_{\bG}}(\wt{\bY})\succeq 0$, $M_k(G_j\wt{\bY})\succeq 0$, $j\in[s]$, and $\rank(M_k(\bY))=\rank(M_{k+d_{\bG}}(\wt{\bY}))$.
\end{enumerate}
\end{theorem}

\subsection{Semialgebraic geometry}
Now, we recall some notions and results of semialgebraic geometry, which can be found in~\cite{realAG}.

\begin{definition}{
A subset $\mathcal{S}$ of $\mathbb{R}^n$ is called semialgebraic if it is a finite union of sets of the form
$$\{\bx \in \mathbb{R}^n \mid f_i(\bx) = 0, \ i = 1, \ldots, p; f_i(\bx) > 0, \ i = p + 1, 
\ldots, q\},$$
where all $f_{i}$ are polynomials.
In other words, $\mathcal{S}$ is a union of finitely many sets, each defined by finitely many polynomial equalities and inequalities.

A mapping $\mathscr{A} \colon \mathcal{S} \rightarrow \mathbb{R}^m$ is said to be semialgebraic if its graph
\begin{equation*}
\{ (\bx, \by) \in \mathcal{S} \times \mathbb{R}^m \mid \by =\mathscr{A}(\bx) \}
\end{equation*}
is a semialgebraic subset of $\mathbb{R}^n \times \mathbb{R}^m.$
A mapping $\mathscr{A}:\mathcal{S}\to\bS^m$ is called semialgebraic
after identifying $\bS^m$ with
$\mathbb{R}^{m(m+1)/2}$.
}\end{definition}

A major fact concerning the class of semialgebraic sets is its stability under linear projections.

\begin{theorem}[\textnormal{Tarski--Seidenberg theorem}] \label{TarskiSeidenbergTheorem}
The image of any semialgebraic set $\mathcal{S} \subset \mathbb{R}^n$ under a projection to any linear subspace of $\mathbb{R}^n$ is a semialgebraic set.
\end{theorem}
As an immediate consequence of the Tarski--Seidenberg Theorem, we get semialgebraicity of any set 
$\{\bx \in \mathcal{A} \mid \exists \by \in \mathcal{B},(\bx, \by) \in \mathcal{C} \},$ provided that $\mathcal{A},\mathcal{B}$ and $\mathcal{C}$ are semialgebraic sets in the corresponding spaces. Also, $\{ \bx \in \mathcal{A} \mid  \forall \by \in \mathcal{B}, (\bx, \by) \in \mathcal{C}\}$ is a semialgebraic set as its complement is the union of the complement of $\mathcal{A}$ and the set $\{ \bx \in \mathcal{A} \mid  \exists \by \in \mathcal{B}, (\bx, \by) \not\in \mathcal{C} \}.$ Thus, if we have a finite collection of semialgebraic sets, then any set obtained from them with the help of a finite chain of quantifiers is also semialgebraic.

\begin{remark}\label{rem::semialg_integrability}
Every finite-valued semialgebraic mapping is Borel measurable. In
particular, if
\[
\mathscr{A}=[a_{ij}]:\K\to\bS^m
\]
is a bounded semialgebraic matrix-valued mapping on $\K$ defined in \eqref{eq::K}, 
then each entry $a_{ij}$ is a bounded Borel measurable function. Consequently, for
every finite signed Borel measure $\phi$ on $\K$,
\[
\int_{\K}|a_{ij}(\bx)|\,\ud|\phi|(\bx)
\le
\|a_{ij}\|_{\infty,\K}\,|\phi|(\K)<\infty,
\]
where $|\phi|$ denotes the total variation measure of $\phi$.
Thus, the pairing $\langle\mathscr{A},\Phi\rangle$ is well defined for
every $\Phi\in\mathfrak{M}^m_+(\K)$.
\end{remark}

\section{A correlative sparse matrix Positivstellensatz}\label{sec::csPos}
In our previous work \cite[Example~5.1]{MWG2026}, we gave a counterexample
showing that the matrix analogue of the scalar Positivstellensatz for
correlative sparsity can fail even under the RIP and local Archimedean
assumptions. In this section, we establish this matrix Positivstellensatz
under an additional pointwise commutation condition.
\subsection{Central common lower envelopes}
As a basic technical tool for establishing this Positivstellensatz, we first construct a central common
lower envelope for a nonempty bounded family of symmetric matrices.
This envelope preserves the infimum of the smallest eigenvalue after adding any Hermitian matrix that commutes with every member of the family.

Let $\{A_u,  u\in\mathcal{U}\}\subset\bS^m$ be a nonempty family of real symmetric matrices and
let 
\[
\mathcal{A}=\mbox{alg}^*_{\CC}\{I_m, A_u, u\in\mathcal{U}\}\subset M_m(\CC)
\]
be the complex unital $*$-algebra generated by $I_m$ and $\{A_u,  u\in\mathcal{U}\}$.
Let $\mathcal{A}'$ be the commutant of $\mathcal{A}$, that is,
\[
\mathcal{A}'=\{B\in M_m(\CC) \mid BA=AB\quad\text{for every}\quad A\in \mathcal{A}\},
\]
and $Z(\mathcal{A})=\mathcal{A}\cap\mathcal{A'}$
be the center of $\mathcal{A}$. Recall that for $A, B\in M_m(\CC)$, the commutator of $A$ and $B$ is
\[
[A, B]=AB-BA.
\]

We have the following standard structure theorem for finite-dimensional
complex $*$-algebras; the formula for the center follows
by intersecting the algebra with its commutant.
\begin{theorem}{\upshape (see \cite[Theorem~5]{KLV2000}, \cite[Theorem~2.1]{MM2011})}\label{th::ast}
    There is a unitary matrix $V\in M_m(\CC)$ and integers $d_{\ell}, e_{\ell}\in\N$ with
    $\sum_{\ell}d_{\ell}e_{\ell}=m$ such that 
    \[
    \begin{aligned}
        V^*\mathcal{A}V&=\bigoplus_{\ell}\left(M_{d_{\ell}}(\CC)\otimes I_{e_{\ell}}\right),\\
        V^*\mathcal{A}'V&=\bigoplus_{\ell}\left(I_{d_{\ell}}\otimes M_{e_{\ell}}(\CC)\right),\\
        V^*Z(\mathcal{A})V&=\bigoplus_{\ell}\CC I_{d_{\ell}e_{\ell}}.\\
    \end{aligned}
    \]
\end{theorem}

\begin{lemma}\label{lem::eigen}
    Let $A=A^*\in M_d(\CC)$ and $R=R^*\in M_e(\CC)$, then 
    \[
    \lambda_{\min}(A\otimes I_e+I_d\otimes R)=\lambda_{\min}(A)+\lambda_{\min}(R).
    \]
\end{lemma}
\begin{proof}
    Note that if $A\bu=\gamma \bu$ and $R\bv=\rho \bv$ for some 
    $\gamma, \rho\in\RR$ and 
    $\bu\in\CC^d$, $\bv\in\CC^e$, then it holds
    \[
    (A\otimes I_e+I_d\otimes R)(\bu\otimes \bv)=(\gamma+\rho)\bu\otimes \bv.
    \]
    Taking orthonormal eigenbases of $A$ and $R$, with
eigenvalues $\gamma_i$ and $\rho_j$, respectively,
their Kronecker products form an orthonormal eigenbasis
of $\CC^{de}$.
Thus all eigenvalues of $A\otimes I_e+I_d\otimes R$
are the pairwise sums $\gamma_i+\rho_j$,
and the claimed identity follows.
\end{proof}

Consider the set
\[
\mathcal{L}:=\{L\in\bS^m\ \mid\  [A_u, L]=0,\ A_u\succeq L \ \text{for all } u\in\mathcal{U}\}.
\]
\begin{proposition}\label{prop::clv}
Suppose that $\mathcal{U}\neq\emptyset$ and  $\sup_{u\in\mathcal{U}}\Vert A_u\Vert\le C<\infty$. Then, 
there is a unique matrix $N$ satisfying
\begin{equation}\label{eq::ccle}
  N\in\mathcal{L},\quad  N\succeq L\quad\text{for all}\quad 
 L\in\mathcal{L}.
\end{equation}
Moreover, $N\in Z(\mathcal A)$, and the following statements hold.
\begin{enumerate}[\upshape (i)]
    \item $-C I_m\preceq N\preceq C I_m$ and $A_u\succeq N$ for all $u\in\mathcal{U}$;
    \item for every Hermitian $R\in\mathcal{A}'$, it holds that $[N, R]=0$ and 
    \[
    \inf_{u\in\mathcal{U}}\lambda_{\min}(A_u+R)=\lambda_{\min}(N+R).
    \]
\end{enumerate}
\end{proposition}
\begin{proof}
It is clear that if a matrix $N$ satisfies \eqref{eq::ccle},
it must be unique. Now let us construct such a matrix.

        Let $V$ be the unitary matrix in Theorem \ref{th::ast}. 
    Then for any $u\in\mathcal{U}$, there exists $A_{\ell}(u)\in M_{d_\ell}(\CC)$ such that 
    \[
    V^*A_uV=\bigoplus_{\ell}\left(A_{\ell}(u)\otimes I_{e_{\ell}}\right).
    \]
    Let $\mu_{\ell}=\inf_{u\in\mathcal{U}}\lambda_{\min} A_{\ell}(u)$.
    As $\sup_{u\in\mathcal{U}}\Vert A_u\Vert\le C$, we have
$-C\le \lambda_{\min}(A_\ell(u))\le C$.
Since $\mathcal U\neq\emptyset$, the set
$\{\lambda_{\min}(A_\ell(u)):u\in\mathcal U\}$
is nonempty and bounded. Therefore, $\mu_\ell$ is well defined and belongs to $[-C,C]$ for all $\ell$. 

Consider a larger set
\[
\wt{\mathcal{L}}
=
\left\{
L\in M_m(\CC)\ \mid\ L^*=L,\ 
[A_u,L]=0,\ A_u-L\succeq 0,\ \forall u\in\mathcal U
\right\}.
\]
For any $L\in\wt{\mathcal{L}}$, we have $L\in\mathcal{A}'$. 
Hence, by Theorem \ref{th::ast}, it holds that
\[
V^*LV
=\bigoplus_\ell\bigl(I_{d_\ell}\otimes L_\ell\bigr)
\quad\text{for some } L_\ell\in M_{e_\ell}(\CC)\quad \text{with}\ L_\ell=L_\ell^*.
\]
Since $A_u\succeq L$ for all $u\in\mathcal{U}$, 
by Lemma \ref{lem::eigen}, we have
\[
\lambda_{\min}A_\ell(u)\ge\lambda_{\max}L_\ell
\quad\forall\,u,\ell.
\]
Therefore, $L_\ell\preceq \mu_\ell I_{e_\ell}$ for all $\ell$.  
Let $N$ be the Hermitian matrix such that 
    \[
    V^*NV=\bigoplus_{\ell}\mu_{\ell}I_{d_{\ell}e_{\ell}}.
    \]
Clearly, $N\in\wt{\mathcal{L}}$ and $N\succeq L$ for all $L\in\wt{\mathcal{L}}$. 

To show \eqref{eq::ccle},
it suffices to prove that $N$ is real. 
Because all $A_u$, $u\in\mathcal{U}$, are real, entrywise complex conjugation preserves $\wt{\mathcal{L}}$. In fact, 
for any $L\in\wt{\mathcal{L}}$, we 
have $\overline{L}^*=\overline{L}$ and 
\[
    [A_u, \overline{L}]=A_u\overline{L}-\overline{L}A_u=\overline{\overline{A_u}L-L\overline{A_u}}=\overline{A_uL-LA_u}=0,\quad 
    A_u=\overline{A_u}\succeq\overline{L},\quad
    \forall\ u\in\mathcal{U}.
\]
Consequently, we have $\overline{N}\in\wt{\mathcal{L}}$. This means that $N\succeq \overline{N}$. Taking entrywise complex conjugation of this inequality yields $\overline{N}\succeq N$, which implies $N=\overline{N}$. Hence, $N$ is real and \eqref{eq::ccle} holds for $N$.

By Theorem \ref{th::ast} and the construction of $N$, 
we also have $N\in Z(\mathcal{A})$.

(i) It is clear from the construction of $N$.

(ii) For any Hermitian $R\in\mathcal{A}'$, as $N\in Z(\mathcal{A})\subset\mathcal{A}$, it is clear that $[N,R]=0$. By Theorem \ref{th::ast}, write $V^*RV=\bigoplus_{\ell}(I_{d_{\ell}}\otimes R_{\ell})$
    with $R_\ell\in M_{e_\ell}(\CC)$.
    As $R$ is Hermitian, each $R_{\ell}$ is also Hermitian.
    For each $u\in\mathcal{U}$, as $A_u$ is symmetric, each $A_{\ell}(u)$ is Hermitian. 
    Thus, it holds
        \begin{align*}
        \lambda_{\min}(A_u+R)&=\lambda_{\min}(V^*A_uV+V^*RV)\\
    &=\min_{\ell}\{\lambda_{\min}(A_{\ell}(u)\otimes I_{e_{\ell}}+
    I_{d_{\ell}}\otimes R_{\ell})\}\\
    &=\min_{\ell}\{\lambda_{\min}A_{\ell}(u)+\lambda_{\min}R_{\ell}\}
    \tag*{[by Lemma~\ref{lem::eigen}]}.
    \end{align*}
    It follows that 
    \begin{align*}
        \inf_{u\in\mathcal{U}}\lambda_{\min}(A_u+R)&=\inf_{u\in\mathcal{U}}\min_{\ell}\{\lambda_{\min}A_{\ell}(u)+\lambda_{\min}R_{\ell}\}\\
        &=\min_{\ell}\left\{\inf_{u\in\mathcal{U}}\lambda_{\min}A_{\ell}(u)+\lambda_{\min}R_{\ell}\right\}\\
        &=\min_{\ell}\left\{\mu_{\ell}+\lambda_{\min}R_{\ell}\right\}\\
        &=\min_{\ell}\lambda_{\min}(\mu_{\ell}I_{d_{\ell}e_{\ell}}+
    I_{d_{\ell}}\otimes R_{\ell})\tag*{[by Lemma~\ref{lem::eigen}]}\\
    &=\lambda_{\min}(V^*NV+V^*RV)=\lambda_{\min}(N+R).
    \end{align*}
The proof completes.
\end{proof}
We call the unique matrix $N$ satisfying \eqref{eq::ccle} the {\itshape central common lower envelope}
of $\{A_u, u\in\mathcal{U}\}$.

%\begin{lemma}\label{lem::length}
%Let $\mathcal B\subseteq M_m(\RR)$ be the unital real
%algebra generated by an arbitrary family $\mathcal G$.
%Then $\mathcal B$ is spanned by words of length at most $m^2-1$
%in the generators, where the word of length zero is $I_m$.
%Moreover, every element of $\mathcal B$ is a real linear
%combination of at most $m^2$ such words.
%\end{lemma}
%\begin{proof}
%For $k\ge0$, let $V_k$ be the real span of all words of length
%at most $k$. Then
%\[
%V_0\subseteq V_1\subseteq\cdots\subseteq\mathcal B,
%\qquad
%\dim V_0=1,\qquad \dim V_k\le m^2.
%\]
%If $V_k=V_{k+1}$, then $V_k$ is invariant under right
%multiplication by every generator. Since $I_m\in V_k$,
%it follows that every word belongs to $V_k$, and hence
%$V_k=\mathcal B$.
%
%Each strict inclusion increases the dimension by at least one.
%Therefore, $V_k=V_{k+1}$ for some $k\le m^2-1$, which gives
%$\mathcal B=V_{m^2-1}$.
%Finally, a basis of $\mathcal B$ can be selected from these
%spanning words and contains at most $m^2$ elements.
%\end{proof}

\begin{proposition}\label{prop::semi}
    Let $K\subset\RR^{n_1+n_2}$ be a nonempty semialgebraic set,
    $\pi_{\RR^{n_1}}(K)$ be the projection of $K$ on the first $n_1$ coordinates, and
    \[
    \mathcal{U}_s=\{u\in\RR^{n_2}\mid (s, u)\in K\}\quad
    \text{for}\quad s\in \pi_{\RR^{n_1}}(K).
    \]
    Let $A: K\to\bS^m$ be a bounded and semialgebraic map. Define a map $M: \pi_{\RR^{n_1}}(K)\to \bS^m$, $s\mapsto N_s$, where $N_s$ denotes the unique central common lower envelope of 
$\{A(s,u), u\in\mathcal{U}_s\}$. 
Then, the map $M$ is bounded, semialgebraic and Borel measurable.
\end{proposition}
\begin{proof}
Since $K$ is semialgebraic and nonempty, $\pi_{\RR^{n_1}}(K)$ is semialgebraic and 
$\mathcal U_s$ is nonempty for every $s\in \pi_{\RR^{n_1}}(K)$.

For $s\in \pi_{\RR^{n_1}}(K)$, define
\[
\mathcal L_s
=
\left\{
L\in\bS^m\ \mid\ 
[A(s,u), L]=0,\ 
A(s,u)\succeq L
\ \text{for all }u\in\mathcal U_s
\right\}.
\]
Symmetry, commutation, and positive semidefiniteness are semialgebraic conditions.
Since the map $A$ is semialgebraic, we conclude that $\{(s,L):s\in \pi_{\RR^{n_1}}(K),\ L\in\mathcal L_s\}$ 
is semialgebraic.

By assumption, choose $C<\infty$ such that $\|A(s,u)\|\le C$ on $K$.
Proposition~\ref{prop::clv} gives, for every $s\in \pi_{\RR^{n_1}}(K)$,
a unique greatest element $N_s$ of $\mathcal L_s$, satisfying
\[
-CI_m\preceq N_s\preceq CI_m.
\]
Therefore, the map $M$ is bounded. By the uniqueness of $N_s$, the graph of $M$ is exactly
\[
\left\{
(s,X)\in \pi_{\RR^{n_1}}(K)\times\bS^m:
X\in\mathcal L_s,\ 
\forall L\in\bS^m\
\bigl(L\in\mathcal L_s\Longrightarrow X-L\succeq0\bigr)
\right\}.
\]
By Tarski--Seidenberg theorem again, $M$ is a semialgebraic map.
It also implies that each entry of the map $M$ is semialgebraic and Borel.
As $M$ is bounded, $M$ is Borel measurable by Remark~\ref{rem::semialg_integrability}.
\end{proof}

%Note that the map $M: S\to \bS^m$ may not continuous.

\subsection{An inequality for consistent matrix-valued measures under commutation}
We next use the preceding envelope construction to derive an integral inequality for local matrix-valued measures with consistent separator marginals.
The resulting inequality bounds the sum of the local integrals from below by the infimum of the smallest eigenvalue of the sum of the local matrix functions over the global feasible set.
No global extension of the local measures is required.

\begin{definition}\label{def::rip}
For a clique family $\{\mI_i\}_{i=1}^t$ of a graph
$\mG=([n],\mE)$ with $\bigcup_{i=1}^t\mI_i=[n]$, we say
    $\{\mI_i\}_{i=1}^t$ satisfies the running intersection property (RIP) if there exists a tree $\mathcal{T}([t], E_{\mathcal{T}})$ rooted at some $i_0\in[t]$ such that for every $j\in[n]$,
the nodes $\{i\in[t]\mid j\in\mI_i\}$ induce a connected subtree.
\end{definition}

\begin{remark}
Besides the clique-tree formulation of the RIP in
Definition~\ref{def::rip}, there is an equivalent ordered formulation.
After a suitable relabeling, for every $\ell=2,\ldots,t$, there exists
an index $p(\ell)<\ell$ such that
\[
\mathcal I_\ell\cap
\left(\bigcup_{j<\ell}\mathcal I_j\right)
\subseteq \mathcal I_{p(\ell)}.
\]
To obtain this formulation from the clique-tree formulation, choose a root and order every
parent before its children. Children at the same level may
be ordered arbitrarily. 
Conversely, choose one admissible index $p(\ell)$ for each
$\ell=2,\ldots,t$, without requiring uniqueness, and connect $\ell$ to
$p(\ell)$. 
\end{remark}

For a clique $\mI_i=\{i_1<\cdots<i_r\}\subseteq[n]$ and
$\bx=(x_1,\ldots,x_n)\in\RR^n$, we write
\[
\bx_{\mI_i}:=(x_{i_1},\ldots,x_{i_r})\in\RR^{\mI_i}
\]
for the subvector of $\bx$ indexed by $\mI_i$. The ordering of the
coordinates is always inherited from the natural ordering of $[n]$.
For $S\subseteq \mI_i$, let
\[
\pi_S:\mathbb{R}^{\mI_i}\to\mathbb{R}^{S},
\qquad
\pi_S(\bx_{\mI_i})=\bx_S
\]
denote the coordinate projection onto the variables indexed by $S$.
If $\Phi_i=[\phi_{pq}]$ is a matrix-valued measure on a set 
$\K_i\subseteq\mathbb{R}^{\mI_i}$, its pushforward under $\pi_S$ is the
matrix-valued measure $(\pi_S)_\#\Phi_i$ on $\mathbb{R}^{S}$ defined by
\[
\bigl((\pi_S)_\#\Phi_i\bigr)(B)
=
\Phi_i\bigl(\K_i\cap\pi_S^{-1}(B)\bigr),
\qquad B\in \mathscr{B}(\mathbb{R}^{S}).
\]
Equivalently,
\[
(\pi_S)_\#\Phi_i
=
\bigl[(\pi_S)_\#\phi_{pq}\bigr].
\]

For a tree $\mathcal{T}=([t],E_{\mathcal{T}})$
rooted at $i_0$, let $\p(i)$ denote the parent of $i$
for $i\neq i_0$, and let $\ch(i)$ denote the set of
children of $i$ for every $i\in[t]$.

\begin{theorem}\label{th::ineq}
Let $\{\mI_i\}_{i=1}^t$ be a clique family of a graph
$\mG=([n],\mE)$ with $\bigcup_i\mathcal I_i=[n]$.
For each $i\in[t]$, let
$\K_i\subseteq\RR^{\mI_i}$,
$A_i:\K_i\to\bS^m$, and
$\Phi_i\in\mathfrak M_+^m(\K_i)$.
Set
\[
\K:=
\{\bx\in\RR^n \ \mid\
  \bx_{\mI_i}\in\K_i,\ i\in[t]\}.
\]
Suppose that the following conditions hold
for $(\mI_i, \K_i, A_i, \Phi_i)_{i=1}^t$:
\begin{enumerate}[\upshape (a).]
    \item The family $\{\mI_i\}_{i=1}^t$ satisfies the RIP
    with respect to a fixed tree
    $\mathcal T=([t],E_{\mathcal T})$ rooted at $i_0$;

    \item The sets $\K$ and $\K_i$, $i\in[t]$,
    are nonempty, compact, and semialgebraic;

    \item Each $A_i$ is bounded and semialgebraic, and
    \[
    [A_i(\bx_{\mI_i}),A_j(\bx_{\mI_j})]=0
    \qquad
    \text{for all }\bx\in\K
    \text{ and }i,j\in[t],\ i\neq j;
    \]

    \item The measures satisfy
    \begin{align}
    &(\pi_{S_i})_\#\Phi_i
      =(\pi_{S_i})_\#\Phi_{\p(i)}
      \qquad
      \text{for all }i\neq i_0
      \text{ with }S_i\neq\emptyset,
      \label{eq::cc}\\
    &\Phi_i(\K_i)=\Phi_{i_0}(\K_{i_0})
      \quad(i\in[t]),
      \qquad
      \tr{\Phi_{i_0}(\K_{i_0})}=1,
      \label{eq::tm}
    \end{align}
    where $S_i:=\mI_i\cap\mI_{\p(i)}$ for $i\neq i_0$.
\end{enumerate}
Then
\begin{equation}\label{eq::main_inequality}
\sum_{i=1}^t\langle A_i,\Phi_i\rangle
\ge
\inf_{\bx\in\K}
\lambda_{\min}\left(
  \sum_{i=1}^t A_i(\bx_{\mI_i})
\right)
=:a^\star.
\end{equation}
\end{theorem}
\begin{proof}
    We prove the inequality by induction on the number $t$ of cliques. 

    \noindent Case $t=1$: We have $A_1(\bx_{\mI_1})-a^{\star}I_m\succeq 0$
    on $\K_1$. Thus, 
    \[
    \langle A_1, \Phi_1\rangle\ge \langle a^{\star} I_m, \Phi_1\rangle=
    a^{\star}\tr{\Phi_1(\K_1)}=a^{\star}.
    \]
    The inequality holds.

    \noindent Case $t>1$: Fix a nonroot leaf node $\ell$ and its parent $p=\p(\ell)$. Let $S=\mI_{\ell}\cap \mI_p$, $U=\mI_{\ell}\setminus S$ and $J=\cup_{j\neq\ell} \mI_j$.

    \vskip 5pt
    We first assume that both $S$ and $U$ are nonempty. After
    possible reordering leaf coordinates, let
    \[
    D_S:=\pi_{S}(\K_{\ell}),\quad \mathcal{U}_s:=\{u \mid (u, s)\in \K_\ell\},\ s\in D_S.
    \]
    Since $A_\ell$ is bounded and semialgebraic,
    applying Proposition~\ref{prop::semi} with $K=\K_{\ell}$ and $A=A_\ell$, we obtain 
    a central lower envelope map $M_{\ell} : D_S\to \bS^m$ which is 
    bounded, semialgebraic and Borel measurable. 
    Next we delete the node $\ell$ and make the following update:
    \begin{equation}\label{eq::update}
        \begin{aligned}
        &\wt{\K}_p=\K_p\cap(\pi_S|_{\RR^{\mI_p}})^{-1}(D_S),\quad
    \wt{A}_p=A_p+M_{\ell}\circ\pi_{S},\quad
    \wt{\Phi}_p=\Phi_p|_{\wt{\K}_{p}},\\
    &\wt{\K}_j=\K_j,\quad \wt{A}_j=A_j,\quad \wt{\Phi}_j=\Phi_j,\quad 
    j\in[t],\ \ j\neq\ell, p,\\
    &\wt{\K}=\{\by\in\RR^{J}\mid \by_{\mI_j}\in\wt{\K}_j,\ j\in[t],\ j\neq \ell\}.
    \end{aligned}
    \end{equation}
    For induction, we now verify that the update \eqref{eq::update} preserves the conditions (a)-(d).

    (a): It is clear that deleting the leaf $\ell$ preserves RIP for the remaining cliques $\{\mI_i, i\in[t], i\neq\ell\}$.

    (b): Since $\K_\ell$ is compact, the projection $D_S$ is compact and
    semialgebraic. Therefore, $\wt{\K}_p$ is a closed subset of the compact $\K_p$ and thus $\wt{\K}_p$ is compact semialgebraic. By the RIP, we have 
    $\mI_\ell\cap J=S$ and $U\cap J=\emptyset$.
    We now prove that $\wt{\K}=\pi_{J}(\K)$. For any $\bx\in\K$, let $\by=\pi_J(\bx)$.
    Then, $\by_{\mI_j}=\bx_{\mI_j}\in\K_j=\wt{\K}_j$ for all $j\neq\ell, p$.
    As $\by_S=\pi_S(\bx_{\mI_\ell})\in\pi_S(\K_\ell)=D_S$, we have $\by_{\mI_p}\in\wt{\K}_p$ 
    (which also shows that $\wt{\K}_p$ is nonempty).
    Then, $\by\in\wt{\K}$ and thus $\pi_{J}(\K)\subseteq\wt{\K}$. Conversely, for any $\by\in\wt{\K}$, 
    as $\by_{\mI_p}\in\wt{\K}_p$, we have $\by_S\in D_S$ and hence there exists $u$
    such that 
    $(u, \by_S)\in\K_{\ell}$. Letting $\bx=(u,\by)$, it implies $\bx_{\mI_j}\in\K_j$ for each $j\in[t]$. Then, we have $\bx\in\K$ and hence $\wt{\K}\subseteq\pi_{J}(\K)$.
    Since $\K$ is nonempty and compact, so is $\wt{\K}$. 

    (c): Since $M_\ell$ is bounded and semialgebraic on $D_S$, by the definition of $\wt{A}_p$
    and $\wt{\K}_p$, the map $\wt{A}_p$ is bounded semialgebraic on $\wt{\K}_p$.
    Fix a point $\by\in\wt{\K}$. Since $\wt{\K}=\pi_{J}(\K)$, we have $\by_S\in D_S$. 
    Then for any $u\in\mathcal{U}_{\by_S}$, 
    we have $(u, \by)\in\K$. 
    Let 
\begin{equation}\label{eq::As}
\mathcal{A}_{\by_S}:=\mbox{alg}^*_{\CC}\{I_m, A_\ell(u,\by_S), u\in\mathcal{U}_{\by_S}\}\subset M_m(\CC),
\end{equation}
and $\mathcal{A}'_{\by_S}$ be its commutant.
    Then for all $j\neq \ell$,
    by assumption, we have $[A_j(\by_{\mI_j}), A_\ell(u, \by_S)]=0$ for all $u\in\mathcal{U}_{\by_S}$
    and hence $A_j(\by_{\mI_j})\in \mathcal{A}'_{\by_S}$. Since $M_\ell(\by_S)\in Z(\mathcal{A}_{\by_S})\subset\mathcal{A}_{\by_S}$, 
    it holds that $[M_\ell(\by_S), A_j(\by_{\mI_j})]=0$ for all $j\neq\ell$. 
    Since every $\by\in\wt{\K}$ extends to a point of $\K$,
the matrices $A_j(\by_{\mI_j})$, $j\neq\ell$,
still commute pairwise.
    Therefore, we have $[\wt{A}_i(\by_{\mI_i}), \wt{A}_j(\by_{\mI_j})]=0$ for all $i, j\in [t]\setminus\{\ell\}$ and $i\neq j$.

    (d): Since $(\pi_{S})_{\#}\Phi_p=(\pi_{S})_{\#}\Phi_\ell$, for the set $V=\K_p\setminus\wt{\K}_p$,
    it holds that
    \[
    \Phi_p(V)=(\pi_{S})_{\#}\Phi_p(\RR^S\setminus D_S)=(\pi_{S})_{\#}\Phi_\ell(\RR^S\setminus D_S)=0_{m\times m}.
    \]
    Hence, for any Borel $E\subset \K_p$, it holds 
    $0\preceq \Phi_p(E\cap V)\preceq \Phi_p(V)=0_{m\times m}$ and hence $\Phi_p(E\cap V)=0_{m\times m}$.
    Then, 
    \[
    \Phi_p(E)=\Phi_p((E\cap V)\cup (E\cap\wt{\K}_p))=\Phi_p(E\cap\wt{\K}_p)=\wt{\Phi}_p(E).
    \]
    Therefore, we have $\wt{\Phi}_p=\Phi_p$. 
    As $\wt{\Phi}_j=\Phi_j$ for all $j\in[t]\setminus\{\ell\}$, $j\neq p$, 
    the consistency condition \eqref{eq::cc} and the total-mass
    condition \eqref{eq::tm} still hold for $(\wt{\Phi}_j)_{j\in[t]\setminus\{\ell\}}$.

    Next, we show that after the update \eqref{eq::update}, the left-hand side of the 
    inequality \eqref{eq::main_inequality} does not increase and the right-hand side remains the same.

    {\itshape Left-hand side:} For any $\bz\in\K_\ell$, by Proposition \ref{prop::semi} and the definition of $M_\ell$, it holds $A_\ell(\bz)\succeq M_\ell(\pi_S(\bz))$. Therefore, we have
    \[
    \begin{aligned}
        \langle A_\ell, \Phi_\ell\rangle&\ge 
    \int_{\K_\ell}\tr{M_\ell(\pi_S(\bz))\ud\Phi_\ell(\bz)}\\
    &=\int_{D_S}\tr{M_\ell(s)\ud((\pi_S)_\#\Phi_\ell)(s)}\\
    &=\int_{D_S}\tr{M_\ell(s)\ud((\pi_S)_\#\Phi_p)(s)}\\
    &=\int_{\wt{\K}_p}\tr{M_\ell(\pi_S(\bz))\ud\wt{\Phi}_p(\bz)}
    =\langle M_\ell\circ\pi_S, \wt{\Phi}_p\rangle.
    \end{aligned}
    \]
    Therefore, as $\wt{\Phi}_j=\Phi_j$ for all $j\in[t], j\neq\ell$, it holds that 
    \begin{equation}\label{eq::lh}
        \begin{aligned}
        \sum_{j\in[t], j\neq\ell} \langle \wt{A}_j, \wt{\Phi}_j\rangle
    &=\langle \wt{A}_p, \wt{\Phi}_p\rangle+
    \sum_{j\in [t], j\neq \ell, p} \langle \wt{A}_j, \wt{\Phi}_j\rangle\\
    &=\langle A_p, \Phi_p\rangle+\langle M_\ell\circ\pi_S, \wt{\Phi}_p\rangle+
    \sum_{j\in [t], j\neq\ell, p} \langle A_j, \Phi_j\rangle\\
    &\le \sum_{j=1}^t  \langle A_j, \Phi_j\rangle.
    \end{aligned}
    \end{equation}

    {\itshape Right-hand side:} 
    Fix a point $\by\in\wt{\K}$ and let $R(\by)=\sum_{j\in[t], j\neq\ell}A_j(\by_{\mI_j})$. 
    Since $\wt{\K}=\pi_{J}(\K)$, we have $\by_S\in D_S$. 
    Then for any $u\in\mathcal{U}_{\by_S}$, we have $(u, \by)\in\K$. By assumption,
    $[A_j(\by_{\mI_j}), A_\ell(u,\by_S)]=0$ for all $j\neq\ell$. Hence, 
    $[R(\by), A_\ell(u,\by_S)]=0$, which means $R(\by)\in\mathcal{A}'_{\by_S}$. Then, by Proposition \ref{prop::clv}, 
    \[
    \inf_{u\in\mathcal{U}_{\by_S}}\lambda_{\min}(A_\ell(u,\by_S)+R(\by))
    =\lambda_{\min}(M_\ell(\by_S)+R(\by))
    =\lambda_{\min}\left(\sum_{j\in[t], j\neq\ell}\wt{A}_j(\by_{\mI_j})\right).
    \]
    Consequently, as $\wt{\K}=\pi_{J}(\K)$, we have
    \begin{equation}\label{eq::rh}
        a^{\star}=\inf_{\bx\in\K}\lambda_{\min}\sum_{j=1}^t A_j(\bx_{\mI_j})
    =\inf_{\by\in\wt{\K}}\inf_{u\in\mathcal{U}_{\by_S}}\lambda_{\min}(A_\ell(u,\by_S)+R(\by))=\inf_{\by\in\wt{\K}}\lambda_{\min}\left(\sum_{j\in[t], j\neq\ell}\wt{A}_j(\by_{\mI_j})\right).
    \end{equation}

    Now we can establish the result. Because the conditions (a)-(d) hold for the update \eqref{eq::update},
    by the induction hypothesis, we have 
    \[
    \sum_{j\in[t], j\neq\ell} \langle \wt{A}_j, \wt{\Phi}_j\rangle\ge 
    \inf_{\by\in\wt{\K}}\lambda_{\min}\left(\sum_{j\in[t], j\neq\ell}\wt{A}_j(\by_{\mI_j})\right).
    \]
    Then, by \eqref{eq::lh} and \eqref{eq::rh}, we can conclude the inequality \eqref{eq::main_inequality}.

\vskip 5pt
    Now it remains to consider the case $U=\emptyset$ and the case $S=\emptyset$.
    
    We first assume that $U=\emptyset$. Then, $S=\mI_\ell$
    and $D_S=\K_\ell$. In the update \eqref{eq::update}, set $D_S=\K_\ell$, $M_\ell=A_\ell$, 
    and leave all other definitions unchanged. 
    This update simply merges the leaf into its parent, thus (a)--(d), \eqref{eq::lh} and \eqref{eq::rh} remain valid,
with equality in \eqref{eq::lh}, and the induction argument applies.

    Finally, assume that $S=\emptyset$. By RIP, we have $\mI_\ell\cap J=\emptyset$.
    Let 
    \[
    \K_{[t]\setminus\{\ell\}}:=\{\by\in\RR^J \mid \by_{\mI_j}\in\K_j,\ j\in[t],\ j\neq\ell\}.
    \]
 Then, we have $\K=\K_\ell\times\K_{[t]\setminus\{\ell\}}$.
Applying Proposition \ref{prop::clv} to the family
$\{A_\ell(\bz)\mid \bz\in\K_\ell\}$ yields a constant matrix
$M_\ell\in\bS^m$. For each $\by\in\K_{[t]\setminus\{\ell\}}$, 
because $\K=\K_\ell\times\K_{[t]\setminus\{\ell\}}$, by assumption, 
it holds 
\[
\left[A_j(\by_{\mI_j}), A_\ell(\bz)\right]=0,\quad \text{for all}\ \bz\in\K_\ell, j\neq\ell.
\]
Therefore, Proposition \ref{prop::clv} yields $\left[A_j(\by_{\mI_j}), M_\ell\right]=0$ for all
$j\neq\ell$ and 
\[
\inf_{\bz\in\K_\ell}
\lambda_{\min}\left(
 A_\ell(\bz)+\sum_{j\neq\ell}A_j(\by_{\mI_j})
\right)
=
\lambda_{\min}\left(
 M_\ell+\sum_{j\neq\ell}A_j(\by_{\mI_j})
\right).
\]
Consequently, 
\begin{equation}\label{eq::1}
    a^{\star}=
\inf_{\by\in\K_{[t]\setminus\{\ell\}}}\inf_{\bz\in\K_\ell}
\lambda_{\min}\left(
 A_\ell(\bz)+\sum_{j\neq\ell}A_j(\by_{\mI_j})
\right)
=
\inf_{\by\in\K_{[t]\setminus\{\ell\}}}\lambda_{\min}\left(
 M_\ell+\sum_{j\neq\ell}A_j(\by_{\mI_j})
\right).
\end{equation}
On the other hand, the lower-bound property of the constant 
$M_\ell$ and the condition \eqref{eq::tm} imply
\begin{equation}\label{eq::2}
    \langle A_\ell,\Phi_\ell\rangle\ge
\langle M_\ell,\Phi_\ell\rangle
=\tr{M_{\ell}\int_{\K_\ell}\ud\Phi_\ell}
=\tr{M_{\ell}\int_{\K_p}\ud\Phi_p}
=\langle M_\ell,\Phi_p\rangle.
\end{equation}
Then, we have 
\begin{align*}
\sum_{i=1}^t\langle A_i,\Phi_i\rangle
&\ge\langle A_p+M_\ell,\Phi_p\rangle+
\sum_{j\in[t], j\neq \ell, p}\langle A_j,\Phi_j\rangle \tag*{[due to \eqref{eq::2}]}\\
&\ge \inf_{\by\in\K_{[t]\setminus\{\ell\}}}\lambda_{\min}\left(
 M_\ell+\sum_{j\neq\ell}A_j(\by_{\mI_j})
\right)\tag*{[by the induction hypothesis]}\\
&=a^{\star}\tag*{[due to \eqref{eq::1}]},
\end{align*}
where the second inequality is due to the fact that
the remaining family, with $A_p$ replaced by $A_p+M_\ell$ and all local sets and measures unchanged,
satisfies conditions (a)--(d).
This completes the proof.
\end{proof}

\subsection{A matrix Positivstellensatz with correlative sparsity}
As our main result, we next establish a correlative sparse matrix Positivstellensatz 
providing a certificate for positive definiteness of a matrix $F(\bx)\in\bS^m[\bx]$
over the set 
\[
    \K=\{\bx \in \RR^n \mid G_j(\bx) \succeq 0, \quad \forall j \in [s]\}. 
\]
where $G_j\in\bS^{q_j}[\bx],\ j\in[s]$.
Let $\bG:=\{G_j,\ j\in[s]\}$.

We define the {\em correlative sparsity pattern graph} associated with $(F, \bG)$, 
denoted by $\mG_{F,\bG}$, such that the vertex set of $\mG_{F, \bG}$ is $[n]$ and for distinct 
$i,j\in[n]$, $\{i,j\}$ is an edge of $\mG_{F,\bG}$ if and only if at least one of following 
conditions holds:
\begin{enumerate}
    \item[(i)] There exists $\ba\in\supp(F)$ such that $i,j\in\supp(\ba)\coloneqq\{k\in[n]\mid\alpha_k\ne0\}$;
    \item[(ii)] There exists $k\in[s]$ such that $x_i,x_j\in\text{var}(G_k)$, 
    where $\text{var}(G_k)$ is the subset of variables effectively involved in $G_k$.
\end{enumerate}

Consider the correlative sparse assumption imposed on $(F,\bG)$:
\begin{assumption}\label{assump2}
   There exist cliques $\{\mI_i\}_{i=1}^t$ of $\mG_{F,\bG}$ with $\bigcup_i\mathcal I_i=[n]$ and subsets $\{\mJ_i\}_{i=1}^t$ of $[s]$ such that
   \begin{enumerate}[\upshape (i)]
 %  \item $\cup_{j=1}^t \mI_j=[n]$ and $\cup_{j=1}^t \mJ_j=[s]$;
       \item $F(\bx)\in\sum_{i=1}^t\bS^m[\bx_{\mI_i}]$;
       \item The sets $\mJ_1,\ldots,\mJ_t\subseteq[s]$ satisfy
$\bigcup_{i=1}^t\mJ_i=[s]$, and
$G_k\in\bS^{q_k}[\bx_{\mI_i}]$
for every $i\in[t]$ and $k\in\mJ_i$;
       \item The family $\{\mI_i\}_{i=1}^t$ satisfies the RIP
    with respect to a fixed tree
    $\mathcal T=([t],E_{\mathcal T})$ rooted at $i_0$.
   \end{enumerate}
\end{assumption}

%Let 
%\[
%\mathcal{M}:=\left\{\Phi=(\Phi_i)_{i\in[t]}\ \middle|\ 
%    \Phi_i\in\mathfrak{M}_+^m(\K_i),\ i\in[t],\ \eqref{eq::cc}\ \text{and}\ \eqref{eq::tm}\ \text{hold}
%\right\}
%\]

Let 
\[
\K_i=\{\bz\in\RR^{\mI_i}\mid G_j(\bz)\succeq 0\,\ j\in \mJ_i\},\ i\in[t].
\]
Clearly, it holds that
\[
\K=\{\bx\in\RR^n \mid \bx_{\mI_i}\in\K_i,\ i\in[t]\}.
\]
For $i\in[t]$, let $\bG_i:=\{G_j, j\in\mJ_i\}$ and denote by \(\QM^m(\bG_i)\) (resp., \(\QM_k^m(\bG_i)\)) the quadratic module (resp., $k$-th truncated quadratic module) generated by $\bG_i$ in \(\bS^m[\bx_{\mI_i}]\).

Recall the concept of algebraic interior point and the convex
separation theorem.
\begin{definition}
    Let $\mathcal{V}$ be a vector space and let $\mathcal{S}\subset\mathcal{V}$ be a convex set. We 
    say that $\bv\in \mathcal{S}$ lies in the algebraic interior of $\mathcal{S}$ if for any 
    straight line $L$ passing through $\bv$, the point $\bv$ lies in 
    the interior of the intersection $\mathcal{S}\cap L$. 
\end{definition}

\begin{theorem}\cite[III, Corollary 1.7]{ACC}\label{th::sep}
    Let $\mathcal{V}$ be a vector space and let $\mathcal{S}_1, \mathcal{S}_2\subset \mathcal{V}$ be non-empty convex sets such that $\mathcal{S}_1\cap \mathcal{S}_2=\emptyset$. Suppose that $\mathcal{S}_1$ has 
    a non-empty algebraic interior. Then, $\mathcal{S}_1$ and $\mathcal{S}_2$ can be separated 
    by an affine hyperplane.
\end{theorem}

\begin{lemma}\label{lem::ai}
    Suppose that the quadratic module $\MQ^m(\bG_i)$ is Archimedean for each $i\in[t]$, then $I_m$ lies in the  algebraic interior of the convex set $\sum_{i=1}^t\MQ^m(\bG_i)$ in the vector space $\sum_{i=1}^t\bS^m[\bx_{\mI_i}]$. 
\end{lemma}
\begin{proof}
    Fix an arbitrary 
    \[
    H(\bx)=\sum_{i=1}^t H_i(\bx_{\mI_i})\in\sum_{i=1}^t\bS^m[\bx_{\mI_i}].
    \]
    For each $i\in[t]$, by the Archimedean condition, $\K_i$ is compact and hence there exists a sufficiently large real number $c_i>0$ 
    such that $c_iI_m\pm H_i(\bx_{\mI_i})\succ 0$ on
    $\K_i$. Then by Theorem \ref{th::psatz},
    we have $c_iI_m\pm H_i(\bx_{\mI_i})\in\MQ^m(\bG_i)$. 
    Then for all $c\ge \sum_{i=1}^t c_i$, we have
    $cI_m\pm H(\bx)\in\sum_{i=1}^t\MQ^m(\bG_i)$, 
    which shows that $I_m$ lies in the  algebraic interior of $\sum_{i=1}^t\MQ^m(\bG_i)$
    in the vector space $\sum_{i=1}^t\bS^m[\bx_{\mI_i}]$.
\end{proof}

In addition to Assumption~\ref{assump2}, we consider the following assumptions.
\begin{assumption}\label{assump::AC}
The following conditions hold.
    \begin{enumerate}[\upshape (i)]
        \item The quadratic module $\MQ^m(\bG_i)$ is Archimedean for each $i\in[t]$;
        \item There is a decomposition $F(\bx)=\sum_{i=1}^t F_i(\bx_{\mI_i})$ where each $F_i\in\bS^m[\bx_{\mI_i}]$ such that the following pointwise commutation
\begin{equation}\label{eq::pc}
    \left[F_i(\bx_{\mI_i}), F_j(\bx_{\mI_j})\right]=0\quad\forall\ \bx\in\K,\ i\neq j\in[t],
\end{equation}
holds.  
    \end{enumerate}
\end{assumption}
Now we can establish our main result.
\begin{theorem}\label{th::csPosi}
    Suppose that $\K\neq\emptyset$, Assumptions \ref{assump2} and \ref{assump::AC} hold.
    If $F(\bx)$ is positive definite on $\K$, then we have
    $F(\bx)\in\sum_{j=1}^t\MQ^m(\bG_j)$.
\end{theorem}
\begin{proof}   
    To the contrary, suppose that $F(\bx)\not\in\sum_{i=1}^t\MQ^m(\bG_i)$. Then, by Lemma \ref{lem::ai} and Theorem \ref{th::sep}, there is a nonzero linear functional $\wt{\mL}: \sum_{i=1}^t\bS^m[\bx_{\mI_i}]\to \RR$ such that 
    \[
    \wt{\mL}(F)\le 0\quad\text{and}\quad \wt{\mL}(H)\ge 0\quad\text{for all}\quad H\in\sum_{i=1}^t\QM^m(\bG_i).
    \]
    We have $\wt{\mL}(I_m)>0$. In fact, if $\wt{\mL}(I_m)=0$, then for any $H\in\sum_{i=1}^t\bS^m[\bx_{\mI_i}]$,
    as is proved in Lemma \ref{lem::ai}, there exists $c>0$ such that 
    $cI_m\pm H\in\sum_{i=1}^t\MQ^m(\bG_i)$.
    Then,
    \[
    0\le \wt{\mL}(cI_m\pm H)=
    c\wt{\mL}(I_m)\pm\wt{\mL}(H)=\pm\wt{\mL}(H),
    \]
    which implies $\wt{\mL}(H)=0$ and hence 
    $\wt{\mL}$ is a zero linear functional, a contradiction.

   Let $\mL=\wt{\mL}/\wt{\mL}(I_m)$ and $\mL_i=\mL|_{\bS^m[\bx_{\mI_i}]}$ for $i\in[t]$. Then $\mL(I_m)=1$, $\mL(F)\le0$, and
$\mL_i$ is nonnegative on $\QM^m(\bG_i)$.
Since $\QM^m(\bG_i)$ is Archimedean,
by Theorem \ref{th::MM}, there exists
$\Phi_i\in\mathfrak M_+^m(\K_i)$ such that
\[
\mL_i(H)=\langle H,\Phi_i\rangle
\qquad
\text{for all}\quad H\in\bS^m[\bx_{\mI_i}].
\]

Let $\{\mI_i\}_{i=1}^t$ be the clique family in 
Assumption~\ref{assump2} satisfying the RIP
with respect to a fixed tree 
$\mathcal T=([t],E_{\mathcal T})$ rooted at $i_0$.
Now we verify the conditions \eqref{eq::cc} and \eqref{eq::tm}. Fix a nonroot node $i$ with $S_i:=\mI_i\cap\mI_{\p(i)}\neq\emptyset$.
For every $H\in\bS^m[\bx_{S_i}]$,
it holds that
\[
\left\langle H,(\pi_{S_i})_\#\Phi_i\right\rangle
=\mL(H)
=\left\langle H,(\pi_{S_i})_\#\Phi_{\p(i)}\right\rangle.
\]
Both the pushforward measures are supported on the set
$\pi_{S_i}(\K_i)\cup\pi_{S_i}(\K_{\p(i)})$, which is compact by the Archimedean condition.
Then by the Stone--Weierstrass theorem and uniqueness of
finite signed measures on a compact set, the condition \eqref{eq::cc} holds.
For every constant matrix $C\in\bS^m$,
we have $\tr{C\Phi_i(\K_i)}=\mL(C)$.
Since the right-hand side is independent of $i\in[t]$,
all total mass matrices coincide.
Taking $C=I_m$ gives
$\tr{\Phi_i(\K_i)}=\mL(I_m)=1$.
Hence \eqref{eq::tm} holds.

Because each $\QM^m(\bG_i)$ is Archimedean, each $\K_i$
is compact, and hence $\K$ is compact.
Since $\K\neq\emptyset$, each $\K_i$ is also nonempty.
The maps $F_i$ are bounded semialgebraic
maps on $\K_i$ and satisfy pointwise commutation. 

Now we have verified that the conditions (a)-(d)
in Theorem \ref{th::ineq} hold for $(\mI_i, \K_i, F_i, \Phi_i)_{i=1}^t$.
Then, by Theorem \ref{th::ineq}, we have
   \[
   0\ge \mL(F)=\sum_{j=1}^t\langle F_j, \Phi_j\rangle
   \ge \inf_{\bx\in\K}\lambda_{\min} F(\bx)>0,
   \]
    a contradiction.
\end{proof}

\begin{example}\cite[Example~5.1]{MWG2026}
\label{ex:mwg-counterexample}
Let 
\[
F(x_1,x_2)=\begin{bmatrix}
2(x_1-1)^2+2\left(x_2-\frac{3}{2}\right)^2+\frac{1}{2}& 3-2x_2\\ 3-2x_2&2(x_1-2)^2+2\left(x_2-\frac{3}{2}\right)^2+\frac{1}{2}
\end{bmatrix},
\]
and $\K=\left\{(x_1, x_2)\in\RR^2 : G_1(x_1):=4-x_1^2\ge0, \, G_2(x_2):=4-x_2^2\ge0\right\}.$
Let $\mI_1=\mJ_1=\{1\}$ and $\mI_2=\mJ_2=\{2\}$.
It is clear that the correlative sparsity pattern in Assumption \ref{assump2} holds and the quadratic 
module $\MQ^m(\bG_i)$ is Archimedean for $i=1, 2$.
It is proved in \cite{MWG2026} that $\inf_{(x_1,x_2)\in\K}\lambda_{\min}(F(x_1,x_2))=\frac{1}{2}$, but
\[
F(x_1,x_2)-\left(\frac{1}{2}-\varepsilon\right)I_2\not\in\sum_{i=1}^2\MQ^m(\bG_i)
\]
for any $\varepsilon\in(0, \frac{1}{2})$.

Now we verify that the pointwise commutation hypothesis \eqref{eq::pc} fails. Let 
\[
F_1(x_1)=\begin{bmatrix}
2(x_1-1)^2&0 \\ 0 &2(x_1-2)^2
\end{bmatrix},\quad
F_2(x_2)=\begin{bmatrix}
2\left(x_2-\frac{3}{2}\right)^2+\frac{1}{2}& 3-2x_2\\ 3-2x_2&2\left(x_2-\frac{3}{2}\right)^2+\frac{1}{2}
\end{bmatrix}.
\]
Then, any decomposition $F(x_1,x_2)=\sum_{j=1}^2 F_j(x_j)$ has the form $F=(F_1+H)+(F_2-H)$.
Since the matrix $H\in\bS^2[x_1]\cap\bS^2[x_2]$, we have $H\in\bS^2$.
For every $H\in\bS^2$, the coefficient of $x_1x_2$ in $[F_1+H,F_2-H]$ is 
\[
\begin{bmatrix}
0& -8\\ 8&0
\end{bmatrix}.
\]
Thus $[F_1+H,F_2-H]\not\equiv 0$.
Since $\K$ has nonempty interior, the pointwise commutation hypothesis \eqref{eq::pc} fails.

Now let us modify $F(x_1,x_2)$ as
\[
\widehat{F}(x_1,x_2)=\begin{bmatrix}
2(x_1-1)^2+2\left(x_2-\frac{3}{2}\right)^2+\frac{1}{2}& \frac{1}{2}\\ \frac{1}{2}&2(x_1-2)^2+2\left(x_2-\frac{3}{2}\right)^2+\frac{1}{2}
\end{bmatrix}.
\]
Decompose $\widehat{F}$ as
$\widehat{F}(x_1,x_2)=\widehat{F}_1(x_1)+\widehat{F}_2(x_2)$, where
\[
\widehat{F}_1(x_1)=\begin{bmatrix}
2(x_1-1)^2& \frac{1}{2}\\ \frac{1}{2}&2(x_1-2)^2
\end{bmatrix},\quad
\widehat{F}_2(x_2)=\left(2\left(x_2-\frac{3}{2}\right)^2+\frac{1}{2}\right)\begin{bmatrix}
1& 0\\ 0&1
\end{bmatrix}.
\]
It is clear that $[\widehat{F}_1(x_1),\widehat{F}_2(x_2)]\equiv 0$. 
Then, all assumptions in Theorem~\ref{th::csPosi} are satisfied.
It is easy to see that $\widehat{F}(x_1,x_2)\succ 0$ on $\K$. 
Direct expansion yields the polynomial identity
\[
\widehat{F}(x_1,x_2)=2\begin{bmatrix}
x_1-1& \frac{1}{4}\\ -\frac{1}{4}&x_1-2
\end{bmatrix}^{\intercal}
\begin{bmatrix}
x_1-1& \frac{1}{4}\\ -\frac{1}{4}&x_1-2
\end{bmatrix}
+\left(2\left(x_2-\frac{3}{2}\right)^2+\frac{3}{8}\right)I_2
\in\MQ^2(\bG_1)+\MQ^2(\bG_2).
\]
\end{example}

\begin{example}\label{ex:commuting-separator}
Let $\boldsymbol x=(x_1,x_2,x_3)$ and take 
\[
\K=\left\{\bx\in\RR^3\ \middle|\ G_1(x_1):=4-x_1^2\ge 0,\ G_2(x_2):=4-x_2^2\ge 0,\ 
G_3(x_3):=\frac{1}{4}-x_3^2\ge 0\right\}.
\]
Set
\[
q(t):=2\left(t-\frac32\right)^2+\frac12,
\quad
T(x_3):=
\begin{pmatrix}
2x_3&1-x_3^2\\
1-x_3^2&-2x_3
\end{pmatrix},
\]
and
\[
F_1(x_1,x_3):=q(x_1)I_2+\left(x_1-\frac32\right)T(x_3),
\quad
F_2(x_2,x_3):=q(x_2)I_2-2\left(x_2-\frac32\right)T(x_3),
\quad
F:=F_1+F_2.
\]

By letting $\mI_1=\mJ_1=\{1,3\}$ and $\mI_2=\mJ_2=\{2,3\}$, Assumption~\ref{assump2} holds.
The quadratic modules $\MQ^2(\bG_i)$ are Archimedean for $i=1, 2$.
Moreover, it is easy to see $[F_1(x_1,x_3),F_2(x_2,x_3)]\equiv 0$.
Therefore, all assumptions in Theorem~\ref{th::csPosi} are satisfied.

We can verify that
\begin{equation}\label{eq::ex}
    F-\frac3{128}I_2
=
2C_1^{\mathsf T}C_1
+2C_2^{\mathsf T}C_2
+\frac58\left(x_3^2+\frac94\right)
 \left(\frac14-x_3^2\right)I_2,
\end{equation}
where
\[
C_1:=\left(x_1-\frac{3}{2}\right)I_2+\frac14T(x_3),
\qquad
C_2:=\left(x_2-\frac{3}{2}\right)I_2-\frac12T(x_3).
\]
Therefore, we have $F(\bx)\succ 0$ on $\K$. 
Meanwhile, \eqref{eq::ex} also shows that 
$F\in\MQ^2(\bG_1)+\MQ^2(\bG_2)$.
\end{example}

\section{Application to Polynomial Matrix Optimization}\label{sec::pmo}
Consider the polynomial matrix optimization problem of
minimizing the smallest eigenvalue of $F(\bx) \in \bS^m[\bx]$ over all $\bx \in \K$:
\begin{equation*}\label{eq:eig_real}
    \lambda^{\star} \coloneqq \sup_{\lambda \in \RR} \lambda \quad\mathrm{s.t.} \quad
     F(\bx) \succeq \lambda I_m,\ \forall \bx \in \K=\{\bx\in\RR^n\mid G_1(\bx)\succeq 0,\ldots, G_s(\bx)\succeq 0\}\neq\emptyset.\tag{\text{PMO}}
\end{equation*}
Suppose the correlative sparsity pattern in Assumption~
\ref{assump2} holds for \eqref{eq:eig_real}. 
In this section, we apply the matrix Positivstellensatz obtained in Theorem~\ref{th::csPosi} to 
\eqref{eq:eig_real} problems with correlative-sparsity to establish the convergence of sparse moment-SOS relaxations.

Let 
\[
d_i=\lceil\deg{G_i}/2\rceil,\ i\in[s],\quad
k_{\min}\coloneqq\max\{\lceil\deg{F}/2\rceil,d_1,\ldots,d_s\}.
\]
For $k\ge k_{\min}$, consider the $k$-th order correlative sparsity adapted SOS relaxations for 
\eqref{eq:eig_real} 
\begin{equation}\label{eq:eig_sos}
  \lambda^{\cssos}_k:=\sup_{\lambda \in \RR} \lambda \quad\mathrm{s.t.} \quad F(\bx) - \lambda I_m \in \sum_{i=1}^t\MQ^m_k(\bG_i).
\end{equation}

Given a sequence $\bY=(Y_{\ba})_{\ba\in\N^n}\subseteq\bS^m$ and a subset $S\subseteq[n]$, 
for $d\in\N$ and $G\in\bS^{q}[\bx_S]$, 
let $M_d(\bY, S)$ (resp. $M_d(G\bY, S)$)
be the moment (resp. localizing) submatrix obtained from $M_d(\bY)$ (resp. $M_d(G\bY)$) 
by retaining only those block rows and columns indexed by $\ba\in \N^n_d$
with $\supp(\ba)\subseteq S$.

Let $\mI:=\{\mI_1,\ldots,\mI_t\}$. For each relaxation order $k$, define
\[
\N_{2k}^{\mI_i}
:=
\left\{
\ba\in\N^n_{2k}
\;\middle|\;
\supp(\ba)\subseteq\mI_i
\right\},\ i\in[t],\quad
\N_{2k}^{\mI}
:=
\bigcup_{i=1}^t
\N_{2k}^{\mI_i}.
\]
The dual $k$-th order correlative sparsity adapted moment relaxation for \eqref{eq:eig_real} is 
given by
\begin{equation}\label{eq::mPMI::csmom}
\left\{
\begin{aligned}
\lambda_k^{\csmom}\coloneqq\inf_{\bY=(Y_{\ba})_{\ba\in\N^\mI_{2k}}}\ &\mL_{\bY}(F)\\
\mathrm{s.t.}\ &M_k(\bY, \mI_i)\succeq0,\quad i\in[t],\\
&M_{k-d_j}(G_j\bY, \mI_i)\succeq0,\quad j\in \mJ_i, \ i\in[t],\\
&\mL_{\bY}(I_m)=1.
\end{aligned}\right.
\end{equation}
In our previous work \cite{MWG2026}, we established sufficient conditions for the finite convergence of the sparse moment--SOS relaxations \eqref{eq:eig_sos} and \eqref{eq::mPMI::csmom}. These conditions rely on extending the local matrix-valued measures associated with an optimal solution of \eqref{eq::mPMI::csmom} to a global matrix-valued measure.

\subsection{Convergence by pointwise commutation}
We now establish asymptotic convergence under a pointwise commutation condition. We also derive a criterion for finite convergence and extraction of a global minimizer under additional local flatness and measure consistency conditions. Neither result requires the local measures to admit a global extension.
\begin{theorem}\label{th::asym_conv}
    We have $\lambda^{\cssos}_k\le \lambda_k^{\csmom}\le \lambda^{\star}$ for all
    $k\ge k_{\min}$. Moreover, suppose that Assumptions~\ref{assump2} and \ref{assump::AC} hold, then 
    $\lambda^{\cssos}_k\nearrow\lambda^{\star}$ and $\lambda^{\csmom}_k\nearrow\lambda^{\star}$
    as $k\to\infty$.
\end{theorem} 
\begin{proof}
%    Fix an arbitrary $\varepsilon>0$. 
    Fix arbitrary $\bu\in\K$ and $\bv\in\RR^m$ with $\Vert\bv\Vert=1$.
%    such that $\lambda^{\star}\le \bv^{\intercal}F(\bu)\bv\le(\lambda^{\star}+\varepsilon)$.
    Define the matrix-valued measure $\Phi=\bv\bv^{\intercal}\delta_{\bu}$.
    For any $k\ge k_{\min}$, let $\bY=(Y_{\ba})_{\ba\in\N^n_{2k}}\subseteq\bS^m$
    be such that $Y_{\ba}=\int_{\K}\bx^{\ba}\ud\Phi(\bx)$ for all $\ba\in\N^n_{2k}$. 
    Then, $\bY$ is feasible to \eqref{eq::mPMI::csmom} and 
    \[
    \lambda_k^{\csmom}\le \mL_{\bY}(F)=\int_\K \tr{F(\bx)\ud\Phi(\bx)}=\bv^{\intercal}F(\bu)\bv.
    \]
    By taking the infimum of $\bv^{\intercal}F(\bu)\bv$ over all such $\bu$ and $\bv$, we obtain $\lambda_k^{\csmom}\le\lambda^{\star}$. By the 
    weak duality, it holds $\lambda^{\cssos}_k\le \lambda_k^{\csmom}\le \lambda^{\star}$ for all
    $k\ge k_{\min}$.

    Now suppose that Assumptions~\ref{assump2} and \ref{assump::AC} hold. 
    For any $\varepsilon>0$, we have $F(\bx)-(\lambda^{\star}-\varepsilon)I_m\succ 0$ on $\K$.
    Then, by Theorem~\ref{th::csPosi}, there exists $k_0\ge k_{\min}$ such that
    \[
    F(\bx)-(\lambda^{\star}-\varepsilon)I_m\in \sum_{i=1}^t\MQ^m_k(\bG_i),
    \]
    for all $k\ge k_0$. It implies that $\lambda_k^{\cssos}\ge \lambda^{\star}-\varepsilon$
    for all $k\ge k_0$. Then letting $\varepsilon\to 0$,
    we have $\lambda^{\cssos}_k\nearrow\lambda^{\star}$ 
    as $k\to\infty$.
    Then, due to that $\lambda^{\cssos}_k\le \lambda_k^{\csmom}\le \lambda^{\star}$ 
    for all
    $k\ge k_{\min}$, it also holds that $\lambda^{\csmom}_k\nearrow\lambda^{\star}$
    as $k\to\infty$.
\end{proof}

%\subsection{Robust lower frame bounds for graph filter banks}
%\label{sec:graph-filter-application}
\begin{example}[Robust lower frame bounds for graph filter banks]\label{ex::gf}
A lower frame bound quantifies how stably a graph signal can be
recovered from the channel outputs \cite{HVG2011}.
In multiscale graph diffusion, band-pass channels can be obtained by
differencing successive low-pass responses \cite{Tong2024}. Motivated by
edge-weight uncertainty in graph filtering \cite{TSB2024}, we consider
the worst-case lower frame bound under joint perturbations of edge weights
and scale parameters.

Consider a filter bank with \(J\ge3\) scales on an undirected graph with
nonnegative edge weights, \(m\) vertices, and \(q\ge1\) uncertain edges.
Set \(n=q+J\), \(\bx=(x_1,\ldots,x_n)\), and \(S=\{1,\ldots,q\}\).
The vector \(\bx_S\in[-1,1]^q\) parametrizes the uncertain edge weights,
and \(x_{q+j}\) is the calibration parameter for scale \(j\in[J]\).
Let \(L_0\) be the nominal graph Laplacian. For each \(i\in S\),
let \(\boldsymbol b_i\in\RR^m\) be the incidence vector of edge \(i\)
and choose a perturbation radius \(\delta_i>0\) at most its nominal weight.
The perturbed Laplacian is
\[
 L(\bx_S)=L_0+\sum_{i=1}^q\delta_i x_i
                  \boldsymbol b_i\boldsymbol b_i^{\mathsf T}.
\]
Let \(L_-\) and \(L_+\) be the Laplacians obtained by setting all
uncertain edge weights to their lower and upper endpoints, respectively.
Then \(0\preceq L_-\preceq L(\bx_S)\preceq L_+\). Choose a normalization
constant \(\lambda_0>0\) with \(\lambda_0\ge\lambda_{\max}(L_+)\).
The shift \(A(\bx_S)=I_m-L(\bx_S)/\lambda_0\) satisfies
\(0\preceq A(\bx_S)\preceq I_m\).

Let \(\tau(\bx_S)=(\sum_{i=1}^q\delta_i x_i)/(\sum_{i=1}^q\delta_i)\)
be the normalized change in total weighted degree. For each
\(j\in[J]\), choose \(\beta_j\in\RR\) and \(\rho_j>0\) with
\(|\beta_j|+\rho_j\le1\). Define the joint feasible set
\[
 \K=\left\{\bx\in[-1,1]^n\ \middle|\
       (x_{q+j}-\beta_j\tau(\bx_S))^2\le\rho_j^2,\quad j\in[J]\right\}.
\]
For every \(\bx_S\in[-1,1]^q\), the centers
\(x_{q+j}=\beta_j\tau(\bx_S)\) give a feasible choice. Setting
\(\beta_j=0\) and \(\rho_j=1\) for all \(j\) gives the independent box case.

Let \(z\) denote the spectral variable. For \(j\in[J]\), let
\(p_j\in\RR[z,x_{q+j}]\) be the low-pass response at scale \(j\)
and set \(P_j=p_j(A(\bx_S),x_{q+j})\).
For example, \(p_j(z,x_{q+j})\) can be a convex combination of \(z^{d_j}\)
and \(z^{d_j+1}\) for an integer \(d_j\ge1\), with weights that are affine
in \(x_{q+j}\).
Set \(p_0=1\) and \(p_{J+1}=0\), so that \(P_0=I_m\) and \(P_{J+1}=0\).
Define \(H_j=P_j-P_{j+1}\) for \(j=0,\ldots,J\) and let
\[
 F(\bx)=\sum_{j=0}^J H_j(\bx)^{\mathsf T}H_j(\bx),\qquad
 \lambda^\star=\min_{\bx\in\K}\lambda_{\min}(F(\bx)).
\]
Thus \(\lambda^\star\) is the largest constant for which
\(\sum_{j=0}^J\|H_j(\bx)\bv\|^2\ge\lambda^\star\|\bv\|^2\)
holds for every \(\bx\in\K\) and every graph signal \(\bv\in\RR^m\).

At each realization \(\bx\), the channel terms \(H_j^{\mathsf T}H_j\)
are polynomials in the same symmetric shift \(A(\bx_S)\), whose
eigenspaces may vary with the edge weights. The channel terms therefore
commute pairwise. For \(j=1,\ldots,J-1\), the interior channel \(H_j\)
depends only on variables indexed by \(\mI_j=S\cup\{q+j,q+j+1\}\).
The blocks satisfy the RIP with respect to the path
\(\mI_1-\cdots-\mI_{J-1}\).
The endpoint channels \(H_0\) and \(H_J\) use variables from the first
and last blocks, respectively. Set
\[
F_1 = H_0^{\mathsf T}H_0+H_1^{\mathsf T}H_1,\quad
F_j = H_j^{\mathsf T}H_j,\quad 2\le j\le J-2,\quad
F_{J-1} = H_{J-1}^{\mathsf T}H_{J-1}+H_J^{\mathsf T}H_J.
\]
Then \(F=\sum_{j=1}^{J-1}F_j\), where each \(F_j\) depends only on
the variables in \(\mI_j\) and the summands commute pointwise.
Assign each box and calibration constraint to every block
containing all of its variables. 
The box constraints make each local quadratic module
Archimedean. If a block requires additional edges in the correlative
sparsity graph, include the redundant constraint
\(\lvert\mI_j\rvert-\sum_{i\in\mI_j}x_i^2\ge0\) for that block to
make it a clique.
Theorem~\ref{th::asym_conv} then guarantees that the lower bounds
$\{\lambda^{\cssos}_k\}_k$ and $\{\lambda_k^{\csmom}\}_k$
converge to \(\lambda^\star\).
At a fixed admissible order \(k\), the largest sparse and dense matrix
moment blocks have sizes \(m\binom{q+2+k}{k}\) and
\(m\binom{q+J+k}{k}\), respectively. For fixed \(m,q,k\), the largest
sparse block is independent of \(J\), and the number of local moment
blocks grows linearly with \(J\).

We use a three-vertex path with two uncertain edge weights
\(7/8+x_i/8\), \(i=1,2\), and set \(\lambda_0=3\). For \(j\in[J]\), we choose
\[
 p_j(z,x_{2+j})=(1-c_j)z+c_jz^2,\quad
 c_j=\frac{7}{32}+\frac{10}{32}\frac{j-1}{J-1}+\frac{x_{2+j}}{8},\quad
 \beta_j=\frac{j-1}{J-1}-\frac12,\quad \rho_j=\frac14.
\]
Thus \(\tau(\bx_S)=(x_1+x_2)/2\). Let \(\bar{\bx}\) be the point
given by \(x_1=x_2=1\) and \(x_{2+j}=3\beta_j/2\) for \(j\in[J]\).
This point is feasible because \(\tau=1\),
\(\lvert x_{2+j}\rvert\le3/4\), and the calibration deviations
are \(\lvert\beta_j\rvert/2\le1/4\). At this point, the two edges
have unit weight, and the unit vector
\(\bv=(1,0,-1)^{\mathsf T}/\sqrt{2}\) is an eigenvector of \(A\)
with eigenvalue \(2/3\). Since
\(c_j=1/8+(j-1)/[2(J-1)]\) at \(\bar{\bx}\), the channel gains along
\(\bv\) are \(13/36\) for \(H_0\), \(1/[9(J-1)]\) for each interior
channel, and \(19/36\) for \(H_J\). Hence
\[
 \lambda^\star(J)\le U_J
 :=\bv^{\mathsf T}F(\bar{\bx})\bv
 =\frac{265}{648}+\frac{1}{81(J-1)}.
\]
We implemented both of the dense and sparse moment relaxations in MATLAB 
with YALMIP \cite{lofberg2004yalmip} and solved them with MOSEK \cite{mosek92}. 
All computations were performed and timed on a desktop 
with 8-Core Intel i7 3.8 GHz CPUs and 8 GB RAM. 
The first admissible relaxation order is \(k=3\) because \(F\)
has degree six.  Table~\ref{tab:graph-filter-timing} reports the computational
results. 
All computed SDP objectives are within \(6.3\times10^{-8}\) of the
feasible upper bound \(U_J\). For \(J=3\), increasing \(k\) from \(3\) to
\(4\) changes each computed objective value by less than
\(10^{-9}\). At \(k=3\), the largest sparse PSD block has size \(105\)
throughout, and the largest dense block grows from \(168\) to \(360\).
The dense-to-sparse wall-time ratios for \(J=3,4,5\) at this order are
\(2.57\), \(9.49\), and \(33.69\), respectively.

\begin{table}[t]
 \centering
 \caption{Computational results for Example \ref{ex::gf}.}
 \label{tab:graph-filter-timing}
 \setlength{\tabcolsep}{4pt}
 \begin{tabular}{ccccllrr}
 \toprule
 & & \multicolumn{2}{c}{Largest PSD block}
   & \multicolumn{2}{c}{Abs.\ diff.\ from \(U_J\)}
   & \multicolumn{2}{c}{Time(s)} \\
 \cmidrule(lr){3-4}\cmidrule(lr){5-6}\cmidrule(lr){7-8}
 \(J\) & \(k\) & Sparse & Dense & Sparse & Dense & Sparse & Dense \\
 \midrule
 $3$ & $3$ & $105$ & $168$ & \(1.86\times10^{-9}\)  & \(1.41\times10^{-12}\) & $2.237$ & $5.757$ \\
 $3$ & $4$ & $210$ & $378$ & \(1.46\times10^{-9}\)  & \(5.18\times10^{-10}\) & $19.357$ & $106.648$ \\
 $4$ & $3$ & $105$ & $252$ & \(1.29\times10^{-8}\)  & \(6.04\times10^{-9}\)  & $3.531$ & $33.523$ \\
 $5$ & $3$ & $105$ & $360$ & \(2.85\times10^{-10}\) & \(6.22\times10^{-8}\)  & $4.663$ & $157.093$ \\
 \bottomrule
 \end{tabular}
\end{table}
\end{example}

Next we consider the finite convergence of the sparse moment relaxation \eqref{eq::mPMI::csmom}.
For each $i\in[t]$, let $d_{\K_i}\coloneqq\max\{1,\ \lceil\deg{G_j}/2\rceil,\ j\in \mJ_i\}$.
\begin{theorem}\label{th::finite_conv}
Let $\bY=(Y_{\ba})_{\ba\in\N^\mI_{2k}}$ be an optimal solution  
of \eqref{eq::mPMI::csmom} with fixed order $k\ge \max\{1, k_{\min}\}$. 
Assume that for each $i\in[t]$, one has
\begin{equation}\label{eq::FE}
       \rank\ {M_k(\bY, \mI_i)}=\rank\ {M_{k-d_{\K_i}}(\bY, \mI_i)},
\end{equation}
then each $(Y_{\ba})_{\ba\in\N^{\mI_i}_{2k}}$, $i\in[t]$, admits a representing 
matrix-valued measure
\[
\Phi_i=\sum_{a=1}^{N_i}W_{ia}\delta_{\bx^{(i,a)}},
\quad W_{ia}\in\bS^m_+,\quad \bx^{(i,a)}\in\K_i,
\quad a\in[N_i],\ N_i\in\N.
\]

Under Assumptions~\ref{assump2} and \ref{assump::AC}, suppose that for each child-parent pair $(i,\p(i))$ in the tree $\mathcal{T}([t],E_{\mathcal{T}})$ with 
$S_i:=\mI_i\cap\mI_{\p(i)}\neq\emptyset$, it holds that
\begin{equation}\label{eq::atomic_consis}
    \sum_{\{a:\,\pi_{S_i}(\bx^{(i,a)})=s\}}W_{ia}
=
\sum_{\{b:\,\pi_{S_i}(\bx^{(\p(i),b)})=s\}}W_{\p(i)b}
\quad\text{for every }\quad s\in D_{S_i},
\end{equation}
where
\[
D_{S_i}:=
\{\pi_{S_i}(\bx^{(i,a)}) \mid a\in[N_i]\}
\cup
\{\pi_{S_i}(\bx^{(\p(i),b)}) \mid b\in[N_{\p(i)}]\},
\]
and an empty sum denotes the zero matrix.
Then $\lambda_k^{\csmom}=\lambda^{\star}$. 
Moreover, let 
\[
\X_i:=\{\bx^{(i,a)}, a\in[N_i]\},\quad i\in[t],
\quad\text{and}\quad
\X:=\{\bx\in\RR^n \mid \bx_{\mI_i}\in \X_i,\ i\in[t]\}.
\]
Then $\emptyset\neq\X\subset\K$ and there exists at least one $\bx^{\star}\in\X$ such that
$\lambda_{\min}F(\bx^{\star})=\lambda^{\star}$.
\end{theorem}
\begin{proof}
    By \cite[Theorem 5]{GW2025} (see also \cite[Theorem 6.2]{KT2022}), the feasibility of $\bY$
    to \eqref{eq::mPMI::csmom}, together with the rank 
    condition \eqref{eq::FE}, implies the existence of the 
    representing measure $\Phi_i$.

Now suppose that Assumptions~\ref{assump2} and \ref{assump::AC} hold.
    The atomic consistency in \eqref{eq::atomic_consis} 
    ensures the validity of the condition \eqref{eq::cc}.
    Note that $\Phi_i(\K_i)=Y_{\bf 0}$ for all $i\in[t]$ and $\tr{\Phi_i(\K_i)}=\mL_{\bY}(I_m)=1$. Hence,
    the condition \eqref{eq::tm} holds for $\Phi_i$, $i\in[t]$.
    By assumption, the conditions (a)-(d)
in Theorem \ref{th::ineq} hold for $(\mI_i, \K_i, F_i, \Phi_i)_{i=1}^t$. By Theorem \ref{th::ineq}, we have
\begin{equation}\label{eq::finite_conv}
    \sum_{i=1}^t\langle F_i,\Phi_i\rangle
\ge
\inf_{\bx\in\K}
\lambda_{\min}F(\bx)=\lambda^{\star}.
\end{equation}
Then, by Theorem~\ref{th::asym_conv}, to derive 
$\lambda^{\csmom}_k=\lambda^{\star}$, it suffices to show that
$\mL_{\bY}(F)=\sum_{i=1}^t\langle F_i,\Phi_i\rangle$. 
Write
\[
\Gamma:=\bigcup_{i=1}^t\supp(F_i),
\qquad
F_i=\sum_{\ba\in\Gamma}F_{i,\ba}\bx^\ba.
\]
For each $\ba\in\Gamma$, the set
\[
J_\ba:=\{i\in[t]\mid\supp(\ba)\subseteq\mI_i\}
\]
is nonempty and, by RIP, induces a connected subtree of
$\mathcal T$. Choose an $i\in J_\ba$ and define
\[
Z_\ba:=\int_{\K_i}\bx^\ba\ud\Phi_i(\bx).
\]
If \(\ba=\mathbf{0}\), all local integrals equal \(Y_{\mathbf{0}}\).
If \(\ba\neq\mathbf{0}\), every separator on the path between any
two nodes of \(J_\ba\) contains \(\supp(\ba)\), so \eqref{eq::atomic_consis}
gives equality of the local integrals along the path.
Thus \(Z_\ba\) is well defined.
Note that
$Z_\ba=Y_\ba$ whenever $\ba\in\Gamma$ and $|\ba|\le2k$.
Moreover, since $F=\sum_iF_i$ and $\deg F\le2k$,
\[
\sum_{i=1}^tF_{i,\ba}=F_\ba=0
\qquad\text{whenever }|\ba|>2k.
\]
Consequently,
\begin{equation}\label{eq::LF}
   \begin{aligned}
\sum_{i=1}^t\langle F_i,\Phi_i\rangle
=\sum_{\ba\in\Gamma}
\tr{\left(\sum_{i=1}^tF_{i,\ba}\right)Z_\ba}
=\sum_{\substack{\ba\in\Gamma\\|\ba|\le2k}}
\tr{F_\ba Y_\ba}
=\mL_{\bY}(F).
\end{aligned} 
\end{equation}
Together with \eqref{eq::finite_conv} and the bound
$\lambda_k^{\csmom}\le\lambda^\star$ from
Theorem~\ref{th::asym_conv}, this proves
$\lambda_k^{\csmom}=\lambda^\star$.

Consider the sets $\X_i$, $i\in[t]$, and $\X$.
We first show that $\X\neq\emptyset$.
For the root node $i_0$ of $\mathcal{T}([t],E_{\mathcal{T}})$, since $\tr{\Phi_{i_0}(\K_{i_0})}=1$,
we can select an atom $\bx^{(i_0, a_{i_0})}$ with nonzero weight $W_{i_0a_{i_0}}$. Now fix a nonroot node $i$, and suppose that an atom
$\bx^{(\p(i),a_{\p(i)})}$ has already been selected in its
parent, with $W_{\p(i)\,a_{\p(i)}}\neq0$.
If $S_i\neq\emptyset$, set
$s_i:=\pi_{S_i}\bigl(\bx^{(\p(i),a_{\p(i)})}\bigr)$.
By \eqref{eq::atomic_consis}, we have
\[
\sum_{\{a:\,\pi_{S_i}(\bx^{(i,a)})=s_i\}}W_{ia}
=
\sum_{\{b:\,\pi_{S_i}(\bx^{(\p(i),b)})=s_i\}}W_{\p(i)b}
\succeq W_{\p(i)\,a_{\p(i)}}\neq0.
\]
Hence we can select an atom $\bx^{(i, a_i)}$ with
$W_{ia_i}\neq0$ and $\pi_{S_i}\bigl(\bx^{(i,a_i)}\bigr)=s_i$.
If $S_i=\emptyset$, select an arbitrary atom $\bx^{(i, a_i)}$ with $W_{ia_i}\neq0$ 
which is possible because $\tr{\Phi_{i}(\K_i)}=1$.
Continuing in this way, we can select $\bx^{(i, a_i)}$ 
for every $i\in[t]$ with matching coordinates indexed in $S_i$.
Then the RIP ensures that the selected atoms $\bx^{(i, a_i)}$ agree on all shared coordinates. 
Consequently, the selected atoms define a global point
$\bx\in\RR^n$ satisfying
\[
\bx_{\mI_i}=\bx^{(i,a_i)}\in\X_i
\qquad\text{for every }i\in[t].
\]
Thus $\bx\in\X$, proving that $\X\neq\emptyset$.
By the definition of $\X_i$'s and $\X$, we have
$\X_i\subset\K_i$, $i\in[t]$, and $\X\subset\K$.
Hence, the pointwise commutation on $\X$
\[
    \left[F_i(\bx_{\mI_i}), F_j(\bx_{\mI_j})\right]=0\quad\forall\ \bx\in\X,\ i\neq j\in[t],
\]
holds. Then, by assumption, the conditions (a)-(d)
in Theorem \ref{th::ineq} hold for $(\mI_i, \X_i, F_i, \Phi_i)_{i=1}^t$. 
By \eqref{eq::LF} and Theorem \ref{th::ineq}, we have
\[
\lambda^{\star}=\mL_{\bY}(F)=\sum_{i=1}^t\langle F_i,\Phi_i\rangle
\ge
\min_{\bx\in\X}
\lambda_{\min}F(\bx)\ge\lambda^{\star}.
\]
Therefore, there exists at least one $\bx^{\star}\in\X$ such that
$\lambda_{\min}F(\bx^{\star})=\lambda^{\star}$.
\end{proof}

\begin{remark}\label{rk::acc}
The atomic consistency condition \eqref{eq::atomic_consis}
is guaranteed if, for each nonroot node $i$ with
$S_i\neq\emptyset$, there exists an integer $r_i$,
$1\le r_i\le k$, such that
\begin{equation}\label{eq::FE2}
\rank M_{r_i}(\bY,S_i)
=
\rank M_{r_i-1}(\bY,S_i).
\end{equation}
Indeed, by \cite[Corollary~5.66]{KT2022}, applied to the flat
extension \eqref{eq::FE2}, the sequence $(Y_\ba\mid \ba\in\N^n_{2r_i},\
\supp(\ba)\subseteq S_i)$
has a unique matrix-valued representing measure.
Since $r_i\le k$, both
$(\pi_{S_i})_\#\Phi_i$ and
$(\pi_{S_i})_\#\Phi_{\p(i)}$
represent this sequence.
Uniqueness therefore implies $(\pi_{S_i})_\#\Phi_i=(\pi_{S_i})_\#\Phi_{\p(i)}$,
which yields \eqref{eq::atomic_consis}.
\end{remark}

Under the flatness condition~\eqref{eq::FE}, the local support
points can be recovered from multiplication matrices as in our previous work ~\cite{GW2025}.
When the degree-\(k\) evaluation matrix at the distinct support
points has full column rank, the corresponding weights follow
from the linear system in~\cite[Remark~3]{GW2025}.

The following example illustrates how Theorem~\ref{th::finite_conv}
allows us to certify finite convergence and extract a global minimizer
without assuming that the given local measures $\Phi_1,\ldots,\Phi_t$
are the marginals of a single positive semidefinite matrix-valued
measure on $\K$.

\begin{example}
Let 
\[
\begin{aligned}
    \K&=\left\{\bx\in\RR^3\ \middle|\ g_1(x_1):=4-x_1^2\ge 0,\ g_2(x_2):=4-x_2^2\ge 0,\ g_3(x_3):=\frac{1}{4}-x_3^2\ge 0\right\}\\
    &=[-2,2]^2\times[-1/2,1/2],
\end{aligned}
\]
and define
\[
 A(x_3)=\begin{pmatrix}
 1-x_3^2&2x_3\\
 2x_3&x_3^2-1
 \end{pmatrix},\qquad
 C=\begin{pmatrix}2&1\\1&1\end{pmatrix},\qquad P(x_1,x_3)=\left(x_1-\frac14\right)I_2-A(x_3).
\]
Consider
\begin{equation}\label{eq:nonmwg-en-problem}
 \lambda^{\star}=\min_{\bx\in \K}\lambda_{\min}(F(\bx)),
\end{equation}
where $F=F_1+F_2$,
\begin{align*}
 F_1(x_1,x_3)=PCP+\frac{1+g_3(x_3)}2I_2,\qquad
 F_2(x_2,x_3)=
 \left[\bigl(x_2(x_2-1+x_3)\bigr)^2
       +\frac{1+g_3(x_3)}2\right]I_2.
\end{align*}
Let $\mI_1=\mJ_1=\{1,3\}$, $\mI_2=\mJ_2=\{2,3\}$. 
It is easy to see that Assumptions~\ref{assump2} and \ref{assump::AC} hold for this problem.
In particular, as $F_2$ is a scalar matrix, we have 
$[F_1(x_1,x_3),F_2(x_2,x_3)]\equiv0$.

One can check that 
\[
F_1+F_2-I_2
=PCP+\bigl(x_2(x_2-1+x_3)\bigr)^2I_2+g_3(x_3)I_2
\succeq0,\qquad \forall\ \bx\in\K.
\]
Thus, $\lambda^{\star}\ge 1$. Moreover, we have $\bu=(3/2,\,0,\,1/2)\in\K$ 
and $\lambda_{\min}F(\bu)=1$. Therefore, it holds $\lambda^{\star}=1$. 

Solving the sparse moment relaxation \eqref{eq::mPMI::csmom} with order $k=3$ yielded an 
optimal solution $\bY^{\star}=(Y^{\star}_\ba)_{\ba\in\N^{\mI}_6}$ satisfying 
$\mL_{\bY^{\star}}(F)=1.0000$ and the rank condition 
\[
\begin{aligned}
    \rank\ {M_3(\bY^{\star}, \mI_1)}&=\rank\ {M_{3-d_{\K_1}}(\bY^{\star}, \mI_1)}=4,\\
    \rank\ {M_3(\bY^{\star}, \mI_2)}&=\rank\ {M_{3-d_{\K_2}}(\bY^{\star}, \mI_2)}=8.
\end{aligned}
\]
From the matrix ${M_3(\bY^{\star}, \mI_1)}$, we extract an atomic matrix-valued measure 
$\Phi^{\star}_1=\sum_{a=1}^{4}W_{1a}\delta_{\bx^{(1,a)}}$ with 
\[
\bx^{(1,1)}=(-1.0000, -0.5000), \ \bx^{(1,2)}=(1.5000, -0.5000),\ 
\bx^{(1,3)}=(-1.0000, 0.5000),\ \bx^{(1,4)}=(1.5000, 0.5000),
\]
and
{\footnotesize
\[
W_{11}=\begin{pmatrix}
 0.0648&0.1297\\
 0.1297&0.2593
 \end{pmatrix},\ 
 W_{12}=\begin{pmatrix}
 0.0874&-0.0437\\
 -0.0437&0.0219
 \end{pmatrix},\
W_{13}=\begin{pmatrix}
 0.0826&-0.1652\\
 -0.1652&0.3304
 \end{pmatrix},\
 W_{14}=\begin{pmatrix}
 0.1228&0.0614\\
 0.0614&0.0307
 \end{pmatrix}.
\]}
From the matrix ${M_3(\bY^{\star}, \mI_2)}$, we extract an atomic matrix-valued measure with 
$8$ atoms. Because these 8 atoms are numerically identical in pairs, we combine their weights 
and obtain the measure $\Phi^{\star}_2=\sum_{a=1}^4 W_{2a}\delta_{\bx^{(2,a)}}$ with 
\[
\begin{aligned}
    \bx^{(2,1)}=(-0.0000, -0.5000), \ \bx^{(2,2)}=(1.5000, -0.5000),\ 
\bx^{(2,3)}=(0.0000, 0.5000),\ \bx^{(2,4)}=(0.5000, 0.5000),
\end{aligned}
\]
and
{\footnotesize
\[
W_{21}=\begin{pmatrix}
 0.1192  &  0.0813\\
    0.0813 &   0.2409
 \end{pmatrix},\ 
 W_{22}=\begin{pmatrix}
 0.0331 &   0.0046\\
    0.0046 &   0.0402
 \end{pmatrix},\
W_{23}=\begin{pmatrix}
  0.1102 &  -0.0589\\
   -0.0589 &   0.2014
 \end{pmatrix},\
 W_{24}=\begin{pmatrix}
 0.0953 &  -0.0449\\
   -0.0449 &   0.1597
 \end{pmatrix}.
\]}
Then we check \eqref{eq::FE2} with $S=\mI_1\cap \mI_2=\{3\}$ and obtain
\[
\rank M_2(\bY^{\star},S)
=
\rank M_1(\bY^{\star},S)=4.
\]
Thus, by Remark \ref{rk::acc}, the atomic consistency condition \eqref{eq::atomic_consis}
is guaranteed. By Theorem \ref{th::finite_conv}, there exists at least one optimal solution 
of \eqref{eq:nonmwg-en-problem} in the set $\X$ consisting of the $8$ points
\[
\begin{aligned}
    \bx^{(\star,1)}&=(-1.0000,-0.0000  , -0.5000),\ 
&&\bx^{(\star,2)}=(   -1.0000,    1.5000  , -0.5000),\\
\bx^{(\star,3)}&=(    1.5000 ,  -0.0000  , -0.5000),\ 
&&\bx^{(\star,4)}=(    1.5000 ,   1.5000  , -0.5000),\\
\bx^{(\star,5)}&=(   -1.0000 ,  -0.0000  ,  0.5000),\
&&\bx^{(\star,6)}=(   -1.0000 ,   0.5000  ,  0.5000),\\
\bx^{(\star,7)}&=(    1.5000 ,  -0.0000  ,  0.5000),\
&&\bx^{(\star,8)}=(    1.5000 ,   0.5000  ,  0.5000).
\end{aligned}
\]
They correspond to the rational points
\[
\left\{\bx\in\RR^3\ \middle|\ x_1\in\left\{-1,\frac32\right\},\
(x_2,x_3)\in
\left\{\left(0,-\frac12\right),\ \left(\frac32,-\frac12\right),\ \left(0,\frac12\right),\ \left(\frac12,\frac12\right)\right\}\right\},
\]
One can check that all these eight rational points are 
global minimizers of \eqref{eq:nonmwg-en-problem}.

Now we point out that the whole sequence  $\bY^{\star}=(Y^{\star}_\ba)_{\ba\in\N^{\mI}_6}$
does not admit a representing PSD matrix-valued measure. 
To see that, consider the matrix
\[
Q(\bx)=R(\bx)^{\intercal}R(\bx),\quad
R(\bx)=
\begin{pmatrix}
 45x_1+45x_3^2-225/4&-90x_3-x_2(x_3+\tfrac12)\\
 -90x_3+x_2(x_3+\tfrac12)&45x_1-45x_3^2+135/4
 \end{pmatrix}.
\]
Direct calculation shows that $Q(\bx)\in\bS^2[x_1,x_3]_4+\bS^2[x_2,x_3]_4$.
Hence, we can apply $\mL_{\bY^\star}$ to $Q$ and obtain
$\mL_{\bY^{\star}}(Q)=-0.0650<0$. As $Q$ is an SOS polynomial matrix, it indicates that 
$\bY^{\star}$ can not admit a representing PSD matrix-valued measure. 
%To quantify the effect of perturbing the moments, write
%\[
% Q_{45}=\sum_{\ba}Q_{\ba}\bx^{\ba},\qquad
% \sum_{\ba,p,q}\bigl|(Q_{\ba})_{pq}\bigr|=58037.625.
%\]
%An entrywise perturbation of the sparse moments of size at most
%$\varepsilon$ can therefore increase the pairing by at most
%$58037.625\,\varepsilon$. At $\varepsilon=10^{-6}$ the upper bound
%is approximately $-0.0070093861$, so the negative sign persists throughout this
%neighborhood. The threshold from the reported pairing is approximately
%$1.12077\times10^{-6}$. At $\varepsilon=10^{-5}$ the same bound is
%$0.51532924$, and this witness gives no sign guarantee at that radius.
\end{example}

\subsection{Computing commuting decompositions by linear equations}

%\section{Constructing exact or approximate commuting decompositions}
To apply the preceding results, we need to verify that $F$ admits a clique decomposition
$F(\bx)=\sum_{i=1}^t F_i(\bx_{\mI_i})$ whose summands satisfy the commutation conditions
\eqref{eq::pc}. In this subsection we impose the stronger ambient
polynomial identities
\begin{equation}\label{eq::pc2}
    \left[F_i(\bx_{\mI_i}),F_j(\bx_{\mI_j})\right]\equiv 0,
    \quad
    1\le i<j\le t.
\end{equation}
The two conditions are equivalent
when $\K$ has nonempty interior, since a polynomial vanishing on a nonempty
open set vanishes identically. 
In the following, we propose
a sufficient construction of such a decomposition to verify \eqref{eq::pc2}, 
although the individual summands need not be specified
in the sparse moment--SOS relaxations.

Fix a degree bound $d\ge\deg(F)$ for the summands, and denote the
coefficient of $\bx^\ba$ in $F_i$ by $F_{i,\ba}$. Matching coefficients
in $F=\sum_i F_i$ gives linear equations. For each pair $i<j$, the
coefficient equations obtained from \eqref{eq::pc2} are bilinear in
the coefficients of $F_i$ and $F_j$. When the input coefficients are real
algebraic numbers and the degree bound is fixed, feasibility can be
decided by real quantifier elimination~\cite{realAG}. However, the
computation can be expensive. 
Next we seek a commuting decomposition by fixing a symmetric polynomial
matrix $A$ and requiring $[F_i,A]\equiv 0$ for each $i$. These identities
give linear equations in the unknown coefficients of the summands.
The following proposition shows that they imply \eqref{eq::pc2} whenever
$A$ has distinct eigenvalues at some point of $\RR^n$.
\begin{proposition}\label{prop::ss}
Let $A,H_1,H_2\in\bS^m[\bx]$. If $A(\bx^{(0)})$ has distinct eigenvalues
at some $\bx^{(0)}\in\RR^n$ and $[H_i,A]\equiv0$ for $i=1,2$, then
$[H_1,H_2]\equiv0$.
\end{proposition}
\begin{proof}
    Since the eigenvalue function is continuous, 
    there exists a nonempty open neighborhood $\mathcal{O}$ of $\bx^{(0)}$ such that 
    $A(\bx)$ has distinct eigenvalues for all $\bx\in\mathcal{O}$. Fixing an $\bx\in\mathcal{O}$, by the spectral theorem, there exist
    orthogonal matrix $U_{\bx}$ and diagonal matrix
    $D_{\bx}$ with distinct diagonal entries such that
    $A(\bx)=U_{\bx}^{\intercal}D_{\bx}U_{\bx}$. For $i=1, 2$, by assumption,
    it holds
    \[
    0=[H_i(\bx), A(\bx)]=
    H_i(\bx)U_{\bx}^{\intercal}D_{\bx}U_{\bx}-U_{\bx}^{\intercal}D_{\bx}U_{\bx}H_i(\bx).
    \]
    Thus, 
    \[
        0= U_{\bx}(H_i(\bx)U_{\bx}^{\intercal}D_{\bx}U_{\bx}-U_{\bx}^{\intercal}D_{\bx}U_{\bx}H_i(\bx))U_{\bx}^{\intercal}
    =U_{\bx}H_i(\bx)U_{\bx}^{\intercal}D_{\bx}-D_{\bx}U_{\bx}H_i(\bx)U_{\bx}^{\intercal}
    =[U_{\bx}H_i(\bx)U_{\bx}^{\intercal}, D_{\bx}].
    \]
    Then as $D_{\bx}$ is diagonal matrix with distinct diagonal entries, $U_{\bx}H_i(\bx)U_{\bx}^{\intercal}$
    must be diagonal. Therefore,
    \[
    \begin{aligned}
        0&=[U_{\bx}H_1(\bx)U_{\bx}^{\intercal}, U_{\bx}H_2(\bx)U_{\bx}^{\intercal}]\\
        &=U_{\bx}H_1(\bx)U_{\bx}^{\intercal}U_{\bx}H_2(\bx)U_{\bx}^{\intercal}-U_{\bx}H_2(\bx)U_{\bx}^{\intercal}U_{\bx}H_1(\bx)U_{\bx}^{\intercal}\\
        &=U_{\bx}H_1(\bx)H_2(\bx)U_{\bx}^{\intercal}-U_{\bx}H_2(\bx)H_1(\bx)U_{\bx}^{\intercal}\\
        &=U_{\bx}[H_1(\bx), H_2(\bx)]U_{\bx}^{\intercal}.
    \end{aligned}
    \]
    Consequently, we have $[H_1(\bx), H_2(\bx)]=0$ for all $\bx\in\mathcal{O}$. Because $\mathcal{O}$ is nonempty and open, it holds that $[H_1, H_2]\equiv 0$.
\end{proof}

%For each node $i\in[t]$, let $\mathcal{T}_i$ be the subtree consisting of $i$ and all its descendants, and
%let 
%\[
%P_i=\sum_{j\in\mathcal{T}_i}F_j.
%\]
%By \eqref{eq::decomp}, we have 
%\[
%\sum_{j\in\mathcal{T}_i}\widehat{F}_j=P_i+B_i.
%\]
%
%\begin{theorem}
%    Suppose that $F(\bx^{(0)})$ has distinct eigenvalues at some 
%$\bx^{(0)}\in\RR^n$. Then for $B_i\in\bS^m[\bx_{S_i}]$, $i\in[t]$, the pairwise ambient commutation \eqref{eq::pc2} holds for $\widehat{F}_i$ in \eqref{eq::decomp} if and only if
%\begin{equation}\label{eq::Bi}
%    [P_i+B_i, F]\equiv 0,\quad i\in[t].
%\end{equation} 
%\end{theorem}
%\begin{proof}
%    Suppose that the condition \eqref{eq::pc2} holds for $\widehat{F}_i$'s in \eqref{eq::decomp}. Then,
%    $[\widehat{F}_j, F]\equiv 0$ for all $j\in[t]$ and 
%    hence
%    \[
%0\equiv\left[\sum_{j\in\mathcal{T}_i}\widehat{F}_j, F\right]=[P_i+B_i, F],\quad i\in[t].
%    \]
%
%    Conversely, suppose that \eqref{eq::Bi} holds. Note that
%    \[
%    \widehat{F}_i=(P_i+B_i)-\sum_{k\in\ch(i)}(P_i+B_i),\quad i\in[t].
%    \]
%    Then, by \eqref{eq::Bi}, we have
%    \[
%    \left[\widehat{F}_i, F\right]\equiv 0,\quad i\in[t].
%    \]
%    Thus, by assumption and Proposition \ref{prop::ss},
%    the condition \eqref{eq::pc2} holds.
%\end{proof}

Write
\[
    F(\bx)=\sum_{\ba\in\N^n}F_\ba\bx^\ba,
    \qquad F_\ba\in\bS^m,
\]
and, for the chosen degree bound $d$, let
\[
\Lambda_{d,i}=\{\ba\in\N^n_d \mid \supp(\ba)\subseteq\mI_i\},\quad i\in[t], \quad\Lambda_d=\supp(F)\cup\bigcup_{i=1}^t\Lambda_{d,i}.
\]
For the fixed matrix $A\in\bS^m[\bx]$, which may be constant, consider
the following system in the unknown symmetric matrices $F_{i,\ba}$
for $i\in[t]$ and $\ba\in\Lambda_d$.
\begin{equation}\label{eq::exact_cons}
\left\{
\begin{aligned}
    &F_\ba-\sum_{i=1}^t F_{i,\ba}=0,\ \ba\in\Lambda_d,\\
    &F_{i,\ba}=0_{m\times m},\quad i\in[t],\ \ba\in\Lambda_d\setminus\Lambda_{d,i},\\
    &\left[\sum_{\ba\in\Lambda_d}F_{i,\ba}\bx^{\ba}, A(\bx)\right]\equiv 0,\quad i\in[t].
\end{aligned}\right.
\end{equation}

Matching polynomial coefficients in the last line gives linear equations
because $A$ is fixed. 
%Summing these commutator identities shows that
%$[F,A]\equiv0$ is necessary for feasibility at any degree.

\begin{proposition}\label{prop::commuting_linear}
Suppose that $A$ has simple spectrum at some point of $\RR^n$.
Then every solution of \eqref{eq::exact_cons} yields a
decomposition $F=\sum_{i=1}^t F_i$ satisfying \eqref{eq::pc2}, where
\[
F_i=\sum_{\ba\in\Lambda_d}F_{i,\ba}\bx^\ba
\in\bS^m[\bx_{\mI_i}]_d,\qquad i\in[t].
\]
\end{proposition}
\begin{proof}
The first two lines of \eqref{eq::exact_cons} ensure that
$F=\sum_{i=1}^t F_i$ and that $F_i\in\bS^m[\bx_{\mI_i}]_d$ for each $i$.
The last line gives $[F_i,A]\equiv0$ for each $i$. Since $A$ has simple
spectrum at some point, Proposition~\ref{prop::ss} implies that
$[F_i,F_j]\equiv0$ for all $i<j$.
\end{proof}

\begin{corollary}\label{cor::commuting_linear}
Set $A=F$ in \eqref{eq::exact_cons}. Suppose that either $t=2$ or $F$
has simple spectrum at some point of $\RR^n$. Then the solutions of
\eqref{eq::exact_cons} are precisely the coefficient families of the
decompositions
\[
F=\sum_{i=1}^t F_i,\qquad F_i=\sum_{\ba\in\Lambda_d}F_{i,\ba}\bx^\ba\in\bS^m[\bx_{\mI_i}]_d,\quad i\in[t],
\]
that satisfy \eqref{eq::pc2}.
\end{corollary}
\begin{proof}
For any such decomposition, the first two lines of \eqref{eq::exact_cons}
hold by construction. The last line follows from
\[
[F_i,F]=\sum_{j=1}^t[F_i,F_j]\equiv0,\qquad i\in[t].
\]
Conversely, any solution of \eqref{eq::exact_cons} gives a decomposition
of the stated form. If $F$ has simple spectrum at some point,
Proposition~\ref{prop::commuting_linear} gives \eqref{eq::pc2}.
If $t=2$, then $F=F_1+F_2$, so $[F_1,F_2]=[F_1,F]\equiv0$.
\end{proof}

\begin{remark}
    For a general reference matrix $A$, infeasibility of \eqref{eq::exact_cons}
rules out only those decompositions within the chosen degree bound whose
summands commute identically with $A$. A different reference matrix may
therefore yield a decomposition satisfying \eqref{eq::pc2} within the same
bound. When $A=F$, infeasibility rules out every decomposition satisfying
\eqref{eq::pc2} within that bound, since each summand of such a decomposition
commutes identically with $F$. A decomposition with higher-degree summands
may still exist.
\end{remark}

\begin{example}\label{ex::linear_commuting_decomposition}
Let $n=3$, $\mI_1=\{1,3\}$, $\mI_2=\{2,3\}$, and define
\[
 D(x_3)=\begin{pmatrix}x_3&1\\1&-x_3\end{pmatrix},\qquad
 F(\bx)=(1+x_1^2+x_2^2+x_3^2)I_2+(x_1+x_2)D(x_3).
\]
An arbitrary decomposition of $F$ might not satisfy the commutation condition \eqref{eq::pc2}. 
For example, let
\begin{align*}
 \widetilde F_1(x_1,x_3)
 &=x_1^2I_2+x_1D(x_3)
   +(1+x_3^2)\begin{pmatrix}1&0\\0&0\end{pmatrix},\\
 \widetilde F_2(x_2,x_3)
 &=x_2^2I_2+x_2D(x_3)
   +(1+x_3^2)\begin{pmatrix}0&0\\0&1\end{pmatrix}.
\end{align*}
Although $\widetilde F_1+\widetilde F_2=F$, we have
\[
 [\widetilde F_1,\widetilde F_2]
 =(1+x_3^2)(x_1+x_2)
   \begin{pmatrix}0&1\\-1&0\end{pmatrix}\not\equiv0.
\]

We take $d=2$ and $A=F$. Since there are two cliques,
Corollary~\ref{cor::commuting_linear} shows that \eqref{eq::exact_cons}
characterizes all decompositions of degree at most two satisfying
\eqref{eq::pc2}. Here $F(1,0,0)$ also has distinct eigenvalues, $1$ and $3$.
Because every
monomial involving $x_1$ can occur only in $F_1$, and every
monomial involving $x_2$ can occur only in $F_2$, $F_1$ and $F_2$ must have
the form
\begin{equation}\label{eq::linear_example_ansatz}
 \begin{aligned}
 F_1(x_1,x_3)&=x_1^2I_2+x_1D(x_3)+B(x_3),\\
 F_2(x_2,x_3)&=x_2^2I_2+x_2D(x_3)+(1+x_3^2)I_2-B(x_3),
 \end{aligned}
\end{equation}
where $B\in\bS^2[x_3]$ has degree at most two. 

Write
\[
 B(x_3)
 =\begin{pmatrix}
 \beta_{11}(x_3)&\beta_{12}(x_3)\\
 \beta_{12}(x_3)&\beta_{22}(x_3)
 \end{pmatrix},\quad
 \beta_{ij}(x_3)=\sum_{\ell=0}^2\beta_{ij,\ell}x_3^{\ell}.
\]
For \eqref{eq::linear_example_ansatz}, the two commutation
identities in \eqref{eq::exact_cons} become
\[
 [F_1,A]=(x_1+x_2)[B,D]\equiv0,
 \qquad
 [F_2,A]=-(x_1+x_2)[B,D]\equiv0.
\]
As polynomial identities, these are equivalent to $[B,D]\equiv0$.
A direct multiplication gives
\[
 [B,D]=
 \begin{pmatrix}
 0&\beta_{11}-\beta_{22}-2x_3\beta_{12}\\
 -\beta_{11}+\beta_{22}+2x_3\beta_{12}&0
 \end{pmatrix}.
\]
Comparing the coefficients of $1,x_3,x_3^2,x_3^3$ gives four linear
equations in the nine coefficients of $B$,
\begin{equation}\label{eq::linear_example_equations}
 \left\{
 \begin{aligned}
 \beta_{11,0}-\beta_{22,0}&=0,\\
 \beta_{11,1}-\beta_{22,1}-2\beta_{12,0}&=0,\\
 \beta_{11,2}-\beta_{22,2}-2\beta_{12,1}&=0,\\
 \beta_{12,2}&=0.
 \end{aligned}
 \right.
\end{equation}
One solution of \eqref{eq::linear_example_equations} is
\[
 B(x_3)=\begin{pmatrix}
 \frac12+x_3^2&\frac12x_3\\
 \frac12x_3&\frac12
 \end{pmatrix}.
\]
Substitution yields
\[
 F_1(x_1,x_3)
 =\begin{pmatrix}
 x_1^2+x_1x_3+x_3^2+\frac12&x_1+\frac12x_3\\
 x_1+\frac12x_3&x_1^2-x_1x_3+\frac12
 \end{pmatrix},\ \
 F_2(x_2,x_3)=\begin{pmatrix}
 x_2^2+x_2x_3+\frac12&x_2-\frac12x_3\\
 x_2-\frac12x_3&x_2^2-x_2x_3+x_3^2+\frac12
 \end{pmatrix}.
\]
The summands add to $F$, and $[F_1,F_2]=[F_1,F]=(x_1+x_2)[B,D]\equiv0$.
\end{example}

\section{Conclusions}\label{sec::cons}
The correlative-sparse matrix Positivstellensatz can fail even under the
running intersection property and the local Archimedean assumptions.
Under these assumptions, we prove that a polynomial matrix $F$ that is
positive definite on the feasible set $\K$ belongs to the sum of the local
matrix quadratic modules whenever $F$ admits a clique decomposition whose
summands commute pairwise on $\K$. We use this result to establish
asymptotic convergence of the corresponding sparse moment--SOS relaxations
for polynomial matrix optimization. We also derive a criterion for finite
convergence of the moment relaxation and extraction of a global minimizer
based on local flatness and agreement of the separator marginals of the
local representing matrix-valued measures. Both convergence results hold
without assuming that the local measures admit a global extension.
Finally, given a symmetric reference polynomial matrix with simple spectrum
at some point, we provide a system of linear equations whose solutions
yield commuting clique decompositions of $F$.

\smallskip 
{\footnotesize
\noindent \textbf{Declaration of AI-Assisted Technologies}. 
The author used ChatGPT (OpenAI) to help formulate the pointwise
commutation hypothesis, develop parts of the proofs and examples. 
All AI-assisted content was reviewed,
edited, and verified by the author, who takes full responsibility
for the paper.
}

\bibliographystyle{abbrv}
\bibliography{cs_pmo_references}

\end{document}